\documentclass[12pt,oneside,reqno]{amsart}
\usepackage{graphicx}
\usepackage{mathrsfs}
\usepackage{stmaryrd}
\usepackage{amsfonts}
\usepackage{cite}
\usepackage{enumerate,amsmath,amssymb,amsthm}
\newcommand{\dif}{\mathrm{d}}

\newcommand{\be}{\begin{eqnarray}}
\newcommand{\ee}{\end{eqnarray}}
\newcommand{\ce}{\begin{eqnarray*}}
\newcommand{\de}{\end{eqnarray*}}
\newtheorem{theorem}{Theorem}[section]
\newtheorem{lemma}[theorem]{Lemma}
\newtheorem{remark}[theorem]{Remark}
\newtheorem{definition}[theorem]{Definition}
\newtheorem{proposition}[theorem]{Proposition}
\newtheorem{Example}[theorem]{Example}
\newtheorem{corollary}[theorem]{Corollary}
\def\e{\varepsilon}
\def\t{\theta}
\def\a{\alpha}

\def\d{\delta}

\def\g{\gamma}

\def\[{{\Big[}}
\def\]{{\Big]}}
\def\<{{\langle}}
\def\>{{\rangle}}
\def\({{\Big(}}
\def\){{\Big)}}

\def\no{\nonumber}
\def\bt{\begin{theorem}}
\def\et{\end{theorem}}
\def\bl{\begin{lemma}}
\def\el{\end{lemma}}
\def\br{\begin{remark}}
\def\er{\end{remark}}
\def\bx{\begin{Example}}
\def\ex{\end{Example}}
\def\bd{\begin{definition}}
\def\ed{\end{definition}}
\def\bp{\begin{proposition}}
\def\ep{\end{proposition}}
\def\bc{\begin{corollary}}
\def\ec{\end{corollary}}

\def\cD{{\mathcal D}}

\def\cM{{\mathcal M}}

\def\cP{{\mathcal P}}

\def\mE{{\mathbb E}}

\def\mL{{\mathbb L}}

\def\mN{{\mathbb N}}

\def\mP{{\mathbb P}}

\def\mR{{\mathbb R}}

\def\sF{{\mathscr F}}

\def\s{\sigma}

\def\geq{\geqslant}
\def\leq{\leqslant}

\begin{document}

\allowdisplaybreaks

\title{Path independence for the additive functionals of stochastic Volterra equations with singular kernels and H\"older continuous coefficients}

\author{Huijie Qiao$^{1*}$ and Jiang-Lun Wu$^2$}



\thanks{This work was partly supported by NSF of China (No. 12071071).}

\thanks{$*$ Corresponding author}


\date{}

\dedicatory{1. School of Mathematics,
Southeast University\\
Nanjing, Jiangsu 211189,  China\\
hjqiaogean@seu.edu.cn\\
2. Department of Mathematics, Computational Foundry, Swansea University\\
Bay Campus, Swansea SA1 8EN, UK\\
j.l.wu@swansea.ac.uk}

\begin{abstract} 
In this paper, we are concerned with stochastic Volterra equations with singular kernels and H\"older continuous coefficients. We first establish the well-posedness 
of these equations by utilising the Yamada-Watanabe approach. Then, we aim to characterise the path-independence for additive functionals of these equations. The 
main challenge here is that the solutions of stochastic Volterra equations are not semimartingales nor Markov processes, thus the existing techniques for obtaining the 
path-independence of usual, semimartingale type stochastic differential equations are no longer applicable. To overcome this difficulty, we link the concerned stochastic 
Volterra equations to mild formulation of certain parabolic type stochastic partial differential equations, and further apply our previous results on the path-independence 
for stochastic evolution equations to get the desired result. Finally, as an important application, we consider a class of stochastic Volterra equations whose kernels are related with fractional Brownian motions and derive the path-independence of additive functionals for them.
\end{abstract}

\maketitle \rm

\section{Introduction}

Given a filtered probability space $(\Omega,\sF, (\sF_t)_{t\geq0}, \mP)$, let $(B_t)_{t\geq 0}$ be a $\mR$-valued $(\sF_t)_{t\geq0}$-Brownian motion on it. 
We consider the following stochastic integral equation on $\mR$:
\be
X_t=x_0+\int_0^t (t-s)^{-\a}b(s,X_s)\dif s+\int_0^t (t-s)^{-\a}\sigma(s,X_s)\dif B_s, \quad t\geq 0,
\label{Eq1}
\ee
where $x_0\in\mR$, $0<\a<\frac{1}{2}$, and the coefficients $b: \mR_+\times\mR\mapsto\mR$, $\sigma: \mR_+\times\mR\mapsto\mR$ are all Borel measurable. 
Usually stochastic integral equations like Eq.(\ref{Eq1}) are called stochastic Volterra equations (SVEs for short). SVEs arise in diverse research fields, such as 
population dynamics, spread of epidemics, mathematical finance like stochastic volatility models, and so on (\cite{gjr, glr}). The study of SVEs can be dated back 
to the early works of Berger and Mizel \cite{bm1, bm2}. Since then, the analysis of SVEs has been extended in many directions, see e.g. \cite{cd,CQian,du,HMR,ms, 
MSZ, pro, wang, zhang} and references in. Let us review some works related to our consideration here. When $b, \s$ satisfy non-Lipschitz conditions, Wang \cite{wang} 
established the well-posedness of Eq.(\ref{Eq1}). We notice that the non-Lipschitz conditions in \cite{wang} don't contain H\"older conditions. Recently, if $b$ does not 
depend on $x$ and $\s$ is H\"older continuous in $x$, Mytnik and Salisbury \cite{ms} studied the well-posedness of Eq.(\ref{Eq1}). Furthermore, Pr\"omel and 
Scheffels \cite{ps} investigated general Eq.(\ref{Eq1}), that is, in Eq.(\ref{Eq1}) the kernel $(t-s)^{-\a}$ is replaced by a joint measurable function $K(t,s)$. They 
required that $b$ is Lipschitz continuous in $x$, $\s$ is H\"older continuous in $x$ and the kernel $K(t,s)$ is regular in $t, s$, and furthermore showed the sample 
path regularity, the integrability and the semimartingale property of solutions to SVEs.

In this paper, we require that $b,\s$ are both H\"older continuous in $x$, and study the well-posedness of Eq.(\ref{Eq1}). For the existence, since $b,\s$ are only continuous in $x$, the classical methods, such as the Picard iteration and the successive approximation, don't work. We establish the tightness for Eq.(\ref{Eq1}) to prove the existence for solutions to Eq.(\ref{Eq1}). Besides, for Eq.(\ref{Eq1}), we emphasise that the kernel $(t-s)^{-\a}$ is singular at the point $t$. Therefore, one can't apply the 
Gronwall inequality and the It\^o formula to verify the pathwise uniqueness for solutions to Eq.(\ref{Eq1}). Here our proof for the pathwise uniqueness is inspired by the work 
of Mytnik and Salisbury \cite{ms}. That is, we transform Eq.(\ref{Eq1}) to a stochastic partial differential equation (SPDE for short). And the latter is common and solvable. Based on the relation between them and the pathwise uniqueness of this SPDE, the pathwise uniqueness for solutions to Eq.(\ref{Eq1}) is established.

On the other hand, in recent years there has been increasing interest utilising stochastic differential equations (SDEs for short) to model the evolution of 
events and uncertain phenomena with randomness, in particular in economics, finance, biology and medicine, as well as in 
ecology and environmental data science, and so on. These diverse studies promote and demand further 
development of stochastic analysis with applied features. An important and remarkable example is the 
path-independence of the Girsanov transformation, also named as the transformation of drifts (c.f., e.g., \cite{ikeda}), for 
SDEs. This property was originated from the consideration of market efficiency in mathematical 
finance, see \cite{hc} wherein the path-independent property of the Black-Scholes model (c.f. \cite{Stein}) under a risk neutral 
probability measure, in other words, the Girsanov transformation for linear SDEs, was addressed 
in a naive manner. The mathematical formulation and justification of the 
path-independence of the Girsanov transformation for general SDEs were initiated in \cite{twwy,wy}, 
wherein a characterisation theorem with necessary and sufficient conditions for the path-independence of the Girsanov 
transformation of (non-degenerated) SDEs is established. Since then, there are a number of papers devoted to 
investigating this property for various SDEs, see more recent works \cite{LXZ,qw,qw1,RW,RY}, 
and the recent survey article \cite{wangwu} and references therein. 

Notice that the solution of Eq.(\ref{Eq1}) is in general not a semimartingale nor a Markov process, preventing the usage of It\^o calculus or Feynman-Kac type formulas. 
Thus, it is a natural question whether one can consider the path-independent property for Eq.(\ref{Eq1}). Here we give description and characterization for this question. In this paper, we connect Eq.(\ref{Eq1}) with a parabolic type SPDE and utilise our previous result regarding the path-independence for stochastic evolution equations in \cite{qw}.

Finally, let us describe our motivation of this paper. We note the fact that by integral transformation of stochastic integrals with respect to fractional Brownian motions, one can link them to the stochastic integrals with respect to standard  Brownian motions involving singular kernels, thus SDEs driven by fractional Brownian motions can be treated as SVEs involving singular kernels. This motives us in this paper to study the path-independence of certain additive functionals for SVEs with singular kernels, with a class of SVEs whose kernels are related with fractional Brownian motions as an important special subclass. 

The rest of our paper is organised as follows. In Section \ref{mai}, we present our main results of the present paper. And the proofs of all theorems are given in Section \ref{thsveeu}-\ref{thpathinde}. A new result on characterising the strong solutions of stochastic differential equations driven by fractional Brownian motions is derived in Section \ref{app}. 

{\bf Notations:} We  introduce some notations used in the sequel.

Let $C(\mR)$ be the collection of continuous functions on $\mR$ and $C^2(\mR)$ be the class of functions on $\mR$ having continuous second derivatives. $C_c^\infty(\mR)$ denotes the collection of all real-valued $C^\infty$ functions of compact
supports.

Set for $\lambda\in\mR$, 
\ce
&&\|h\|_{\lambda,\infty}:=\sup\limits_{x\in\mR}|h(x)|e^{-\lambda |x|},\\
&&C_{tem}:=\left\{h\in C(\mR), \|h\|_{\lambda,\infty}<\infty ~\mbox{for every}~ \lambda>0\right\},
\de
and we equip $C_{tem}$ with the topology induced by the norms $\|\cdot\|_{\lambda,\infty}$ for $\lambda>0$. Similarly we define 
$$
C_{rap}:=\left\{h\in C(\mR), \|h\|_{\lambda,\infty}<\infty ~\mbox{for every}~ \lambda<0\right\},
$$
and endow it with the topology induced by the norms $\|\cdot\|_{\lambda,\infty}$ for $\lambda<0$. Let $C^k_{tem}$ (respectively $C^k_{rap}$)  be the collection of functions in $C_{tem}$ (respectively in $C_{rap}$)  which are $k$ times continuously differentiable with all the derivatives in $C_{tem}$ (respectively in $C_{rap}$). Let 
$C(\mR_+, C_{tem})$ be the space of all continuous functions on $\mR_+$ taking values in $C_{tem}$. Let $\cM$ be the space of finite measures on $\mR$ endowed with the weak convergence topology.

The following convention will be used throughout the paper: $C$ with or without indices will denote different positive constants whose values may change 
from one place to another.
 
\section{Main results}\label{mai}

In this section, we introduce our main results of the present paper.

\subsection{SVEs with singular kernels}\label{svesk}

In this subsection, we introduce SVEs with singular kernels and the notion of the path-independence for their additive functionals.

We assume:
\begin{enumerate}[(${\bf H}^1_{b}$)]
\item 
For any $t\geq 0$, $x\mapsto b(t,x)$ is decreasing; and there exist two constants $L_b\geq 0$ and $\g_1\in(0,1]$ such that for any $t\geq 0$ and $x_1, x_2, x\in\mR$
\ce
&&|b(t,x_1)-b(t,x_2)|\leq L_b|x_1-x_2|^{\g_1},\\
&&|b(t,x)|\leq L_b(1+|x|).
\de
\end{enumerate}
\begin{enumerate}[(${\bf H}^1_{\sigma}$)]
\item
There exist two constants $L_\s\geq 0$ and $\g_2\in[\frac{1}{2},1]$ such that for any $t\geq 0$ and $x_1, x_2, x\in\mR$
\ce
&&|\sigma(t,x_1)-\sigma(t,x_2)|\leq L_{\s}|x_1-x_2|^{\g_2},\\
&&|\s(t,x)|\leq L_{\s}(1+|x|).
\de
\end{enumerate}

\bt\label{eu}
Assume that (${\bf H}^1_{b}$) and (${\bf H}^1_{\sigma}$) hold. Then there exists a unique strong solution $(X_t)_{t\geq 0}$ to Eq.(\ref{Eq1}) fulfilling  
\be
\mE|X_t|^p<\infty, \quad \forall p>\frac{2}{1-2\a}.
\label{xbou}
\ee
\et

We will give our proof of the above theorem in Section \ref{thsveeu}.

\br
When $b$ is independent of $x$ and $\s$ is independent of $t$, Theorem \ref{eu} is just right Theorem 1.1 in \cite{ms}. Therefore, our result is more general.
\er

Next, we arbitrarily fix $T>0$ and introduce the following additive functional
\be
f_{s,t}&:=&\int_s^t g_1(r, \omega, X_r)\dif r+\int_s^t g_2(r, \omega, X_r)\dif B_r, \quad 0\leq s< t\leq T,
\label{addfun}
\ee
where 
\ce
g_1: [0,T]\times\Omega\times\mR\mapsto\mR, \qquad g_2:  [0,T]\times\Omega\times\mR\mapsto\mR
\de
are progressively measurable, and moreover for almost all $\omega\in\Omega$, both $g_1(t, \omega, x)$ and $g_2(t, \omega, x)$ are continuous in $(t,x)$, so that $f_{s,t}$ is a well-defined semimartingale. 

\bd\label{pathinde1}
The additive functional $f_{s,t}$ is called path-independent with respect to the solution $(X_t)_{t\in[0,T]}$ of Eq.(\ref{Eq1}), if there exists a function $v: [0,T]\times\mR\mapsto\mR$ such that 
\be
f_{s,t}=v(t,X_t)-v(s,X_s).
\label{defi1}
\ee
\ed

\subsection{A SPDE relating to Eq.(\ref{Eq1})}\label{subspde}

In this subsection, we set up the connection of Eq.(\ref{Eq1}) to a SPDE.

Set $\t:=\frac{1}{\a}-2>0$ so that 
$$
\frac{1}{2+\t}=\a.
$$
Define the operator $\Delta_\t$ and its domain as  
\ce
&&\Delta_\t:=\frac{2}{(2+\t)^2}\frac{\partial}{\partial x}\(|x|^{-\t}\frac{\partial}{\partial x}\),\\
&&\cD(\Delta_\t):=\{\phi\in C_{rap}^2: \Delta_\t\phi\in C_{rap}\},\\
&&\cD_{tem}(\Delta_\t):=\{\phi\in C_{tem}^2: \Delta_\t\phi\in C_{tem}\}.
\de
Next, we consider the following evolution equation on $\mR_+\times\mR$
\ce\left\{ \begin{array}{l}
\frac{\partial}{\partial t}u(t,x)=\Delta_\t u(t,x),\\
u(0,x)=u_0(x)\in\cD(\Delta_\t).
\end{array}
\right.
\de
By direct calculation, one obtains that 
\ce
p^\t_t(x):=\left\{ \begin{array}{l}
c_\t\frac{e^{-\frac{|x|^{2+\t}}{2t}}}{t^\a}, \quad t>0, x\in\mR,\\
\d_0(x), \qquad\quad t=0, x\in\mR,
\end{array}
\right.
\de
is the fundamental solution to the above evolution equation, where the normalised constant $c_\t$ is determined as follows 
$$
c_\t:=\left(\int_{\mR}\frac{e^{-\frac{|x|^{2+\t}}{2t}}}{t^\a}\dif x\right)^{-1}=\left(\int_{\mR}e^{-\frac{1}{2}|x|^{2+\t}}\dif x\right)^{-1},
$$
and $\d_0$ stands for the Dirac delta function, that is, $\d_0(x)=0$, for $x\neq 0$, $\d_0(x)=\infty$, for $x=0$ and $\int_{\mR}\d_0(x)\dif x=1$. Note that $p^\t_t(\cdot)$ is a probability density on $\mR$ for each $t\ge0$. Let $\{S_t^\t, t\geq 0\}$ be the semigroup generated by $\Delta_\t$, namely 
$$
(S^\t_t\phi)(x):=\int_{\mR}p^\t_t(x-y)\phi(y)\dif y, \quad \phi\in C_{tem}^2.
$$

Now let us consider the following SPDE 
\be
\dif {\bf X}_t=\left(\Delta_\t {\bf X}_t+\frac{1}{c_\t}b(t,{\bf X}_t(0))\d_0\right)\dif t+\frac{1}{c_\t}\sigma(t,{\bf X}_t(0))\d_0\dif B_t 
\label{spde}
\ee
for a random field ${\bf X}_t={\bf X}_t(x)={\bf X}(t,x,\omega):\mR_+\times\mR\times\Omega\mapsto\mR$. Then the following expression 
\be
{\bf X}_t=S_t^\t {\bf X}_0+\int_0^tS_{t-s}^\t \frac{1}{c_\t}b(s,{\bf X}_s(0))\d_0\dif s+\int_0^tS_{t-s}^\t\frac{1}{c_\t}\sigma(s,{\bf X}_s(0))\d_0\dif B_s
\label{mildsolu}
\ee
gives an explicit mild solution to the equation (\ref{spde}). Moreover, we say that ${\bf X}$ is a weak solution to the equation (\ref{spde}) if for any $t\in [0,T]$ and $\phi\in \cD(\Delta_\t)$ ${\bf X}$ satisfies
\be
\<{\bf X}_t,\phi\>&=&\<{\bf X}_0,\phi\>+\int_0^t\<{\bf X}_s,\Delta_\t\phi\>\dif s+\int_0^t\<\frac{1}{c_\t}b(s,{\bf X}_s(0))\d_0,\phi\>\dif s\no\\
&&+\int_0^t\<\phi,\frac{1}{c_\t}\sigma(s,{\bf X}_s(0))\d_0\dif B_s\>,
\label{weaksolu}
\ee
where $\<{\bf X}_t,\phi_t\>:=\int_{\mR} {\bf X}_t(y)\phi_t(y)\dif y$.

In the following, we describe the relationship between weak solutions and mild solutions and establish the well-posedness of the equation (\ref{spde}).

\bt\label{equimildweak}
If ${\bf X}$ is a weak solution for the equation (\ref{spde}), ${\bf X}$ is a mild solution to the equation (\ref{spde}). Conversely if ${\bf X}$ is a mild solution for the equation (\ref{spde}), then ${\bf X}$ is also a weak solution to the equation (\ref{spde}).
\et

\bt\label{exisuniq}
Assume that ${\bf X}_0\in C_{tem}$. Suppose further that $b, \sigma$ fulfil (${\bf H}^1_{b}$) (${\bf H}^1_{\sigma}$), respectively. Then, there exists a pathwise unique mild solution ${\bf X}\in C(\mR_+, C_{tem})$ to the equation (\ref{spde}).
\et

The proofs of two above theorems are both placed in Section \ref{thspdeeu}.

\br\label{spdeth}
We mention that a special case of Theorem \ref{exisuniq} has appeared in \cite[Theorem 2.5]{ms}.
\er

Next, let us establish a relation between weak solutions of the equation (\ref{Eq1}) and mild solutions of the equation (\ref{spde}). Take ${\bf X}_0=x_0$ in (\ref{mildsolu}), and notice the following 
\ce
\(S_t^\t {\bf X}_0\)(x)&=&\(S_t^\t x_0\)(x)=x_0,\\
\(S_{t-s}^\t\frac{1}{c_\t}b(s,{\bf X}_s(0))\d_0\)(x)&=&\int_{\mR}p_{t-s}^\t(x-y)\frac{1}{c_\t}b(s,{\bf X}_s(0))\d_0(y)\dif y\\
&=&p_{t-s}^\t(x-0)\frac{1}{c_\t}b(s,{\bf X}_s(0)),\\
\(S_{t-s}^\t\frac{1}{c_\t}\sigma(s,{\bf X}_s(0))\d_0\)(x)&=&\int_{\mR}p_{t-s}^\t(x-y)\frac{1}{c_\t}\sigma(s,{\bf X}_s(0))\d_0(y)\dif y\\
&=&p_{t-s}^\t(x-0)\frac{1}{c_\t}\sigma(s,{\bf X}_s(0)), 
\de
where the property of the Dirac delta function is used, and (\ref{mildsolu}) becomes 
\ce
{\bf X}_t(x)=x_0+\int_0^tp_{t-s}^\t(x-0)\frac{1}{c_\t}b(s,{\bf X}_s(0))\dif s+\int_0^tp_{t-s}^\t(x-0)\frac{1}{c_\t}\sigma(s,{\bf X}_s(0))\dif B_s.
\de
In particular for $x=0$, we get 
\ce
{\bf X}_t(0)&=&x_0+\int_0^tp_{t-s}^\t(0-0)\frac{1}{c_\t}b(s,{\bf X}_s(0))\dif s+\int_0^tp_{t-s}^\t(0-0)\frac{1}{c_\t}\sigma(s,{\bf X}_s(0))\dif B_s\\
&=&x_0+\int_0^t (t-s)^{-\a}b(s,{\bf X}_s(0))\dif s+\int_0^t(t-s)^{-\a}\sigma(s,{\bf X}_s(0))\dif B_s.
\de
That is, ${\bf X}_t(0)$ is a weak solution of the equation (\ref{Eq1}). Thus, in order to obtain the path-independent property for the equation (\ref{Eq1}), 
we are going to study the path-independence related to the equation (\ref{spde}).

\subsection{An auxiliary SDE}\label{auxsde}

In this subsection, based on the equation (\ref{spde}), we want to construct an auxiliary SDE and study its path-independent property.

By Theorem \ref{exisuniq}, we conclude that ${\bf X}_t$ belongs to the Banach space $C_{tem}$. Hence, we will consider a weak solution of the equation (\ref{spde}), i.e. for any $\varphi\in C(\mR_+, \cD(\Delta_\t))$ with $r\mapsto\frac{\partial \varphi_r(\cdot)}{\partial r}\in C(\mR_+, C_{rap})$, 
\be
\<{\bf X}_t,\varphi_t\>&=&\<{\bf X}_0,\varphi_0\>+\int_0^t\left<{\bf X}_r,\left(\Delta_\t\varphi_r+\frac{\partial \varphi_r}{\partial r}\right)\right>\dif r+\int_0^t\frac{1}{c_\t}b(r,{\bf X}_r(0))\varphi_r(0)\dif r\no\\
&&+\int_0^t\frac{1}{c_\t}\sigma(r,{\bf X}_r(0))\varphi_r(0)\dif B_r,
\label{weaksoluspde}
\ee
where $\varphi_t(\cdot): \mR\mapsto \mR $. Based on Theorem \ref{equimildweak} and \ref{exisuniq}, there exists a pathwise unique solution $({\bf X}_t)_{t\geq 0}$ to the equation (\ref{weaksoluspde}). On the other hand, for each $\varphi\in C(\mR_+, \cD(\Delta_\t))$, notice that the equation (\ref{weaksoluspde}) is an It\^o type SDE, thus by the Yamada-Watanabe theorem, the pair 
$\<{\bf X}_t,\varphi_t\>$ is a unique strong solution (in the It\^o sense) of the equation (\ref{weaksoluspde}). Now, we define the path-independent property related to the equation (\ref{weaksoluspde}). 

\bd\label{pathinde2}
For the additive functional
\be
F^{\varphi}_{s,t}:=\int_s^t G_1(r, \omega, \<{\bf X}_r,\varphi_r\>)\dif r+\int_s^t G_2(r, \omega, \<{\bf X}_r,\varphi_r\>)\dif B_r, \quad 0\leq s< t\leq T,
\label{addfun1}
\ee
where 
\ce
G_1: [0,T]\times\Omega\times\mR\mapsto\mR, \qquad G_2: [0,T]\times\Omega\times\mR\mapsto\mR
\de
are progressively measurable, if there exists a function $V: [0,T]\times\mR\mapsto\mR$ such that $F^{\varphi}_{s,t}, V$ satisfy
\be
F^{\varphi}_{s,t}=V(t,\<{\bf X}_t,\varphi_t\>)-V(s,\<{\bf X}_s,\varphi_s\>),
\label{defi2}
\ee
we say that $F^{\varphi}_{s,t}$ is path-independent.
\ed

For the equation (\ref{weaksoluspde}) and the additive functional $F^{\varphi}_{s,t}$, we have the following result.

\bt\label{funthe}
Assume that $\sigma(r,{\bf X}_r(0))\varphi_r(0)\not\equiv 0$ for any $r\in[0,T]$, and for $\omega\in\Omega$, $G_1(\cdot,\omega,\cdot)\in C([0,T]\times\mR\mapsto\mR)$, $G_2(\cdot,\omega,\cdot)\in C([0,T]\times\mR\mapsto\mR)$. Then, the additive functional $F^{\varphi}_{s,t}$ is path-independent with $V$ belonging to $C_b^{1,2}([0,T]\times\mR)$ and all the derivatives of $V$ in $(t,z)$ being uniformly continuous, i.e., the equality (\ref{defi2}) holds if and only if $V, G_1, G_2$ satisfy the following  partial differential equations 
\be\left\{\begin{array}{ll}
\partial_r V(r,z)+\frac{1}{2}\partial_z^2 V(r,z)\left(\frac{1}{c_\t}\sigma(r,{\bf X}_r(0))\varphi_r(0)\right)^2\\
+\partial_z V(r,z)\left[\<{\bf X}_r,\left(\Delta_\t\varphi_r+\frac{\partial \varphi_r}{\partial r}\right)\>+\frac{1}{c_\t}b(r,{\bf X}_r(0))\varphi_r(0)\right]=G_1(r, z),\\
\partial_z V(r,z)\frac{1}{c_\t}\sigma(r,{\bf X}_r(0))\varphi_r(0)=G_2(r, z).
\end{array}
\label{eq2}
\right.
\ee
 \et
 
 The proof of the above theorem is placed in Section \ref{thfunaux}.

\subsection{The path-independence for SVEs}\label{pathindeeq1}

In this subsection, we study the path-independence for the equation (\ref{Eq1}).

\bt\label{pathindeth}
Suppose that $\sigma\not\equiv 0$, for $\omega\in\Omega$, $g_1(\cdot,\omega,\cdot)\in C([0,T]\times\mR\mapsto\mR)$, $g_2(\cdot,\omega,\cdot)\in C([0,T]\times\mR\mapsto\mR)$ and $g_1, g_2$ are uniformly bounded. Then, under (${\bf H}^1_{b}$) (${\bf H}^1_{\sigma}$), the additive functional $f_{s,t}$ is path-independent with respect to $(X_t)_{t\in[0,T]}$ with $v$ belonging to $C_b^{1,2}([0,T]\times\mR)$ and all the derivatives of $v$ being uniformly continuous, if and only if
 \be
&&\partial_r v(r,z)+\frac{1}{2}\partial_z^2 v(r,z)\left(\frac{1}{c_\t}\sigma(r,X_r)\right)^2+\partial_z v(r,z)\frac{1}{c_\t}b(r,X_r)=g_1(r, z),\label{equisequ0}\\
&&\partial_z v(r,z)\frac{1}{c_\t}\sigma(r,X_r)=g_2(r, z).\label{equisequ01}
\ee
\et

The proof of the above theorem is placed in Section \ref{thpathinde}.

\br
We mention that in \cite{qw}, if Eq.(1) doesn't contain the state distribution and jumps, the framework is similar to that here. However, comparing \cite[Theorem 3.2]{qw} with Theorem 2.8, we find that here Eq.(\ref{equisequ0}) and Eq.(\ref{equisequ01}) depend on $X_r$ and $\frac{1}{c_\t}$, which is due to the character of SVEs.
\er

\section{Proof of Theorem \ref{eu}}\label{thsveeu} 

In this section, we present our proof to Theorem \ref{eu}. To this end, our strategy is that we first prove Theorem \ref{eu} under the special case of $\g_1=\g_2=1$ (Proposition \ref{cont}), and for the general case with $\gamma_1\in (0,1], \gamma_2\in [1/2,1]$, establish the existence for solutions of Eq.(\ref{Eq1}) in Proposition \ref{cont1}. Moreover, the argument after Remark \ref{spdeth} together with our Theorem \ref{exisuniq} then assures the uniqueness for solutions of Eq.(\ref{Eq1}).

\bp\label{cont}
Assume that a function $\mR\ni t\mapsto g(t)\in\mR$ is continuous in $t$ and $b, \s$ satisfy (${\bf H}^1_{b}$) (${\bf H}^1_{\sigma}$) with $\g_1=\g_2=1$. Then the following equation
\be
U_t=g(t)+\int_0^t (t-s)^{-\a}b(s,U_s)\dif s+\int_0^t (t-s)^{-\a}\sigma(s,U_s)\dif B_s
\label{ug}
\ee
has a unique solution $(U_t)_{t\geq 0}$ with
\be
\mE|U_t|^p<\infty, \quad p>\frac{2}{1-2\a}.
\label{ubou}
\ee
\ep
\begin{proof}
Fix $T>0$ and construct the Picard iterated sequence as follows: For $n\in\mN$
\ce\left\{ \begin{array}{l}
U_t^{(0)}=g(t)\in\mR,\\
U_t^{(n)}=g(t)+\int_0^t (t-s)^{-\a}b(s,U_s^{(n-1)})\dif s+\int_0^t (t-s)^{-\a}\sigma(s,U_s^{(n-1)})\dif B_s, \quad t\in [0,T].
\end{array}
\right.
\de
For $p>\frac{2}{1-2\a}$, it holds that 
\ce
\mE|U_t^{(n)}|^p&\leq&3^{p-1}|g(t)|^p+3^{p-1}\mE\left|\int_0^t (t-s)^{-\a}b(s,U_s^{(n-1)})\dif s\right|^p\\
&&+3^{p-1}\mE\left|\int_0^t (t-s)^{-\a}\sigma(s,U_s^{(n-1)})\dif B_s\right|^p\\
&\leq&3^{p-1}|g(t)|^p+3^{p-1}\mE\left|\int_0^t (t-s)^{-\a}b(s,U_s^{(n-1)})\dif s\right|^p\\
&&+3^{p-1}C_p\mE\left|\int_0^t (t-s)^{-2\a}|\sigma(s,U_s^{(n-1)})|^2\dif s\right|^{p/2}\\
&\leq& 3^{p-1}|g(t)|^p+3^{p-1}L_b^p\mE\left|\int_0^t (t-s)^{-\a}(1+|U_s^{(n-1)}|)\dif s\right|^p\\
&&+2^{p/2}L_{\s}^p3^{p-1}C_p\mE\left|\int_0^t (t-s)^{-2\a}(1+|U_s^{(n-1)}|^2)\dif s\right|^{p/2}\\
&\leq& 3^{p-1}|g(t)|^p+6^{p-1}L_b^p\mE\left|\int_0^t (t-s)^{-\a}\dif s\right|^p+L_{\s}^p6^{p-1}C_p\left|\int_0^t (t-s)^{-2\a}\dif s\right|^{p/2}\\
&&+6^{p-1}L_b^p\left(\int_0^t (t-s)^{-\frac{\a p}{p-1}}\dif s\right)^{(p-1)}\int_0^t\mE|U_s^{(n-1)}|^p\dif s\\
&&+L_{\s}^p6^{p-1}C_p\left(\int_0^t (t-s)^{-\frac{2\a p}{p-2}}\dif s\right)^{(p-2)/2}\int_0^t\mE|U_s^{(n-1)}|^p\dif s\\
&\leq&C_{p,T}+C_{p,T}\int_0^t\mE|U_s^{(n-1)}|^p\dif s,
\de
where we use the Burkholder-Davis-Gundy inequality, (${\bf H}^1_{b}$) (${\bf H}^1_{\sigma}$), the H\"older inequality and the fact that for $0<\a<\frac{1}{2}$,
\ce
&&\int_0^t (t-s)^{-\a}\dif s\leq C_T,\quad \int_0^t (t-s)^{-2\a}\dif s\leq C_T, \\
&&\int_0^t (t-s)^{-\frac{\a p}{p-1}}\dif s\leq C_T,\quad \int_0^t (t-s)^{-\frac{2\a p}{p-2}}\dif s\leq C_T.
\de
By iteration, we have that
\ce
\mE|U_t^{(n)}|^p\leq C_{p,T}e^{C_{p,T}T}+\sup\limits_{t\in[0,T]}|g(t)|^p\frac{(C_{p,T}t)^n}{n!}\leq \left(C_{p,T}+\sup\limits_{t\in[0,T]}|g(t)|^p\right)e^{C_{p,T}T},
\de
and furthermore
\be
\sup\limits_{n}\sup\limits_{t\in[0,T]}\mE|U_t^{(n)}|^p\leq \left(C_{p,T}+\sup\limits_{t\in[0,T]}|g(t)|^p\right)e^{C_{p,T}T}.
\label{xnboun}
\ee

For $n, m\in\mN$, set $Z_t^{n,m}:=U_t^{(n)}-U_t^{(m)}$, and by the Burkholder-Davis-Gundy inequality, (${\bf H}^1_{b}$) (${\bf H}^1_{\sigma}$) with $\g_1=\g_2=1$ and the H\"older inequality it holds that
\ce
\mE|Z_t^{n,m}|^p&\leq&2^{p-1}\mE\left|\int_0^t (t-s)^{-\a}b(s,U_s^{(n-1)})\dif s-\int_0^t (t-s)^{-\a}b(s,U_s^{(m-1)})\dif s\right|^p\\
&&+2^{p-1}\mE\left|\int_0^t (t-s)^{-\a}\sigma(s,U_s^{(n-1)})\dif B_s-\int_0^t (t-s)^{-\a}\sigma(s,U_s^{(m-1)})\dif B_s\right|^p\\
&\leq&2^{p-1}\mE\left(\int_0^t (t-s)^{-\a}|b(s,U_s^{(n-1)})-b(s,U_s^{(m-1)})|\dif s\right)^{p}\\
&&+2^{p-1}C_p\mE\left(\int_0^t (t-s)^{-2\a}|\sigma(s,U_s^{(n-1)})-\sigma(s,U_s^{(m-1)})|^2\dif s\right)^{p/2}\\
&\leq&2^{p-1}L_{b}^p\left(\int_0^t (t-s)^{-\frac{\a p}{p-1}}\dif s\right)^{(p-1)}\int_0^t\mE|U_s^{(n-1)}-U_s^{(m-1)}|^p\dif s\\
&&+2^{p-1}C_pL_{\s}^p\left(\int_0^t (t-s)^{-\frac{2\a p}{p-2}}\dif s\right)^{(p-2)/2}\int_0^t\mE|U_s^{(n-1)}-U_s^{(m-1)}|^p\dif s\\
&\leq&C_{p,T}\int_0^t\mE|Z_s^{n-1,m-1}|^p\dif s.
\de
Moreover, integrating two sides of the above inequality, we obtain that
\ce
\int_0^t\mE|Z_s^{n,m}|^p\dif s\leq C_{p,T}\int_0^t\left(\int_0^s\mE|Z_r^{n-1,m-1}|^p\dif r\right)\dif s\leq C_{p,T}t\int_0^t\frac{1}{s}\left(\int_0^s\mE|Z_r^{n-1,m-1}|^p\dif r\right)\dif s,
\de
and
$$
\frac{1}{t}\int_0^t\mE|Z_s^{n,m}|^p\dif s\leq C_{p,T}\int_0^t\frac{1}{s}\left(\int_0^s\mE|Z_r^{n-1,m-1}|^p\dif r\right)\dif s.
$$
Again set 
$$
h_t^{n,m}:=\frac{1}{t}\int_0^t\mE|Z_s^{n,m}|^p\dif s, 
$$
and it holds that
\ce 
h_t^{n,m}\leq C_{p,T}\int_0^t h_s^{n-1,m-1}\dif s.
\de
Based on (\ref{xnboun}), it is easy to see that
$$
\sup\limits_{n,m}\sup\limits_{t\in[0,T]}h_t^{n,m}\leq C_{p,T}.
$$
Thus, by the Fatou lemma, one can get that 
$$
h_t\leq C_{p,T}\int_0^t h_s\dif s,
$$
where $h_t:=\limsup\limits_{n,m\rightarrow\infty}h_t^{n,m}$, which together with \cite[Lemma 2.1]{wang} yields that $h_t=0$. This means that $\{U^{(n)}, n\in\mN\}$ is a Cauchy sequence in $\mL^p:=L^p([0,T]\times\Omega, \cP_T, \dif t\times\mP; \mR)$, where 
$\cP_T$ denotes the collection of the progressive measurable sets of $[0,T]\times\Omega$. Since $\mL^p$ is complete, there exists a process $U\in\mL^p$ satisfying
\ce
\lim\limits_{n\rightarrow\infty}\mE\int_0^T|U_t^{(n)}-U_t|^p\dif t=0.
\de

Next, we prove that $U$ is a solution of Eq.(\ref{ug}). Indeed, by the Burkholder-Davis-Gundy inequality, (${\bf H}^1_{b}$) (${\bf H}^1_{\sigma}$) with $\g_1=\g_2=1$ and the H\"older inequality it holds that
\ce
&&\mE\int_0^T\left|\int_0^t (t-s)^{-\a}b(s,U_s^{(n)})\dif s-\int_0^t (t-s)^{-\a}b(s,U_s)\dif s\right|^p\dif t\\
&\leq&\int_0^T\mE\left(\int_0^t (t-s)^{-\a}|b(s,U_s^{(n)})-b(s,U_s)|\dif s\right)^{p}\dif t\\
&\leq&L_{b}^p\int_0^T\left(\int_0^t (t-s)^{-\frac{\a p}{p-1}}\dif s\right)^{(p-1)}\left(\int_0^t\mE|U_s^{(n)}-U_s|^p\dif s\right)\dif t\\
&\leq&C_{p,T}\mE\int_0^T|U_t^{(n)}-U_t|^p\dif t,
\de
and
\ce
&&\mE\int_0^T\left|\int_0^t (t-s)^{-\a}\sigma(s,U_s^{(n)})\dif B_s-\int_0^t (t-s)^{-\a}\sigma(s,U_s)\dif B_s\right|^p\dif t\\
&\leq&C_p\int_0^T\mE\left(\int_0^t (t-s)^{-2\a}|\sigma(s,U_s^{(n)})-\sigma(s,U_s)|^2\dif s\right)^{p/2}\dif t\\
&\leq&L_{\s}^pC_p\int_0^T\left(\int_0^t (t-s)^{-\frac{2\a p}{p-2}}\dif s\right)^{(p-2)/2}\left(\int_0^t\mE|U_s^{(n)}-U_s|^p\dif s\right)\dif t\\
&\leq&C_{p,T}\mE\int_0^T|U_t^{(n)}-U_t|^p\dif t,
\de
which implies that
\ce
&&\int_0^t (t-s)^{-\a}b(s,U_s^{(n)})\dif s\rightarrow \int_0^t (t-s)^{-\a}b(s,U_s)\dif s ~\mbox{in}~\mL^p,\\
&&\int_0^t (t-s)^{-\a}\sigma(s,U_s^{(n)})\dif B_s\rightarrow \int_0^t (t-s)^{-\a}\sigma(s,U_s)\dif B_s ~\mbox{in}~\mL^p.
\de
So, $U$ is a solution of Eq.(\ref{ug}). 

By the similar deduction, we also conclude that the solutions of Eq.(\ref{ug}) are pathwise unique.
\end{proof}

\bp\label{cont1}
Assume that $g(t)$ is continuous in $t$ and $b, \s$ satisfy (${\bf H}^1_{b}$) (${\bf H}^1_{\sigma}$). Then Eq.(\ref{ug}) has a solution $(U_t)_{t\geq 0}$ satisfying (\ref{ubou}).
\ep
\begin{proof}
First of all, we fix any $T>0$ and choose two sequences of functions $\{b^m\}_{m\in\mN}, \{\s^m\}_{m\in\mN}$ which satisfy the following conditions: for any $t\in[0,T]$ and $m\in\mN$
\ce
&&|b^m(t,x_1)-b^m(t,x_2)|+|\s^m(t,x_1)-\s^m(t,x_2)|\leq L^m_1|x_1-x_2|, \quad x_1,x_2\in\mR,\\
&&|b^m(t,x)|+|\s^m(t,x)|\leq L_2(1+|x|), \quad x\in\mR,
\de
where $L^m_1, L_2>0$ are two constants, $L^m_1$ is independent of $t$ and $L_2$ is independent of $m, t$, and
\ce
\lim\limits_{m\rightarrow\infty}\sup\limits_{t\in[0,T],x\in\mR}|b^m(t,x)-b(t,x)|=0, \quad \lim\limits_{m\rightarrow\infty}\sup\limits_{t\in[0,T],x\in\mR}|\s^m(t,x)-\s(t,x)|=0.
\de
Then by Proposition \ref{cont} for each $m\geq1$ there exists a process $U^m$ that solves Eq.(\ref{ug}) with $b^m, \s^m$. In the following, we divide into two steps.

{\bf Step 1.} We prove that $\{U^m\}_{m\in\mN}$ is tight in $C([0,T], \mathbb{R})$.

On one side, by the same deduction to that for (\ref{xnboun}), it holds that for any $p>\frac{2}{1-2\a}$,
\be
\sup\limits_{m}\sup\limits_{t\in[0,T]}\mE|U_t^{m}|^p\leq \left(C_{p,T}+\sup\limits_{t\in[0,T]}|g(t)|^p\right)e^{C_{p,T}T}.
\label{umunbou}
\ee
On the other side, for $0\leq t<t+\t\leq T$ and any $p>\frac{2}{1-2\a}$, 
\be
&&\mE|U^m_{t+\t}-U^m_t|^p\no\\
&=&\mE\bigg|\int_0^{t+\t}(t+\t-s)^{-\a}b^m(s,U_s^m)\dif s+\int_0^{t+\t}(t+\t-s)^{-\a}\s^m(s,U_s^m)\dif B_s\no\\
&&\quad-\int_0^{t}(t-s)^{-\a}b^m(s,U_s^m)\dif s-\int_0^{t}(t-s)^{-\a}\s^m(s,U_s^m)\dif B_s\bigg|^p\no\\
&\leq&4^{p-1}\mE\left|\int_{t}^{t+\t} (t+\t-s)^{-\a}b^m(s,U_s^m)\dif s\right|^p\no\\
&&+4^{p-1}\mE\left|\int_0^{t} \left((t+\t-s)^{-\a}-(t-s)^{-\a}\right)b^m(s,U_s^m)\dif s\right|^p\no\\
&&+4^{p-1}\mE\left|\int_{t}^{t+\t} (t+\t-s)^{-\a}\s^m(s,U_s^m)\dif B_s\right|^p\no\\
&&+4^{p-1}\mE\left|\int_0^{t} \left((t+\t-s)^{-\a}-(t-s)^{-\a}\right)\s^m(s,U_s^m)\dif B_s\right|^p\no\\
&=:&I_1+I_2+I_3+I_4.
\label{i1i2i3i4}
\ee

For $I_1+I_3$, by the H\"older inequality, the Burkholder-Davis-Gundy inequality and (\ref{umunbou}), we get that
\be
I_1+I_3&\leq&4^{p-1}\mE\left|\int_{t}^{t+\t}(t+\t-s)^{-\a}|b^m(s,U_s^m)|\dif s\right|^p\no\\
&&+4^{p-1}C_p\mE\left|\int_{t}^{t+\t} (t+\t-s)^{-2\a}|\s^m(s,U_s^m)|^2\dif s\right|^{p/2}\no\\
&\leq&4^{p-1}\left(\int_{t}^{t+\t} (t+\t-s)^{-\frac{\a p}{p-1}}\dif s\right)^{p-1}\int_{t}^{t+\t}\mE|b^m(s,U_s^m)|^p\dif s\no\\
&&+4^{p-1}C_p \left(\int_{t}^{t+\t} (t+\t-s)^{-\frac{2\a p}{p-2}}\dif s\right)^{\frac{p-2}{2}}\int_{t}^{t+\t}\mE|\s^m(s,U_s^m)|^p\dif s\no\\
&\leq&4^{p-1}C_p\t^{(1-\a)p}+4^{p-1}C_p\t^{(\frac{1}{2}-\a)p}.
\label{i1i3}
\ee
For $I_2+I_4$, by the Burkholder-Davis-Gundy inequality and the extended Minkowski inequality in \cite[Corollary 1.32, P.27]{ka}, it holds that
\ce
I_2+I_4&\leq& 4^{p-1}\mE\left|\int_0^{t} \left|(t+\t-s)^{-\a}-(t-s)^{-\a}\right||b^m(s,U_s^m)|\dif s\right|^p\\
&&+4^{p-1}C_p\mE\left|\int_0^{t} \left|(t+\t-s)^{-\a}-(t-s)^{-\a}\right|^2|\s^m(s,U_s^m)|^2\dif s\right|^{p/2}\\
&\leq& 4^{p-1}C_p\mE\left|\int_0^{t} \left|(t+\t-s)^{-\a}-(t-s)^{-\a}\right|(1+|U_s^m|)\dif s\right|^p\\
&&+4^{p-1}C_p\mE\left|\int_0^{t} \left|(t+\t-s)^{-\a}-(t-s)^{-\a}\right|^2(1+|U_s^m|^2)\dif s\right|^{p/2}\\
&\leq& 4^{p-1}\left(\int_0^{t} \left|(t+\t-s)^{-\a}-(t-s)^{-\a}\right|\(1+(\mE|U_s^m|^p)^{1/p}\)\dif s\right)^p\\
&&+4^{p-1}C_p\left(\int_0^{t} \left|(t+\t-s)^{-\a}-(t-s)^{-\a}\right|^2\(1+(\mE|U_s^m|^p)^{2/p}\)\dif s\right)^{p/2}\\
&\leq& 4^{p-1}C_p\left(\int_0^{t} \left|(t+\t-s)^{-\a}-(t-s)^{-\a}\right|\dif s\right)^p\\
&&+4^{p-1}C_p\left(\int_0^{t} \left|(t+\t-s)^{-\a}-(t-s)^{-\a}\right|^2\dif s\right)^{p/2}.
\de
Note that
\ce
0\leq\int_0^{t} \left((t-s)^{-\a}-(t+\t-s)^{-\a}\right)\dif s=\frac{t^{1-\a}}{1-\a}+\frac{\t^{1-\a}}{1-\a}-\frac{(t+\t)^{1-\a}}{1-\a}\leq \frac{\t^{1-\a}}{1-\a},
\de
and 
\ce
&&\int_0^{t} \left((t-s)^{-\a}-(t+\t-s)^{-\a}\right)^2\dif s=\int_0^{t}\left(\frac{(t+\t-s)^{\a}-(t-s)^{\a}}{(t-s)^{\a}(t+\t-s)^{\a}}\right)^2\dif s\\
&\leq&\t^{2\a}\int_0^{t}\frac{1}{(t-s)^{2\a}(t+\t-s)^{2\a}}\dif s=\t^{1-2\a}\int_0^{\frac{t}{\t}}\frac{1}{(v+1)^{2\a}v^{2\a}}\dif v\\
&\leq&\left\{ \begin{array}{l}C_\a\t^{2\a}, ~ \a\in(0,\frac{1}{4}),\\
C_\a\t^{\frac{1}{2}-\a}, \a\in[\frac{1}{4},\frac{1}{2}).
\end{array}
\right.
\de
Thus, one can obtain that
\be
I_2+I_4\leq \left\{ \begin{array}{l}C\t^{(1-\a)p}+C_\a\t^{\a p}, ~ \a\in(0,\frac{1}{4}),\\
C\t^{(1-\a)p}+C_\a\t^{\frac{1-2\a}{4}p}, \a\in[\frac{1}{4},\frac{1}{2}).
\end{array}
\right.
\label{i2i4}
\ee
Combining (\ref{i1i3}) (\ref{i2i4}) with (\ref{i1i2i3i4}), we conclude that
\be
\lim\limits_{\t\downarrow 0}\sup\limits_{m}\sup\limits_{0\leq t<t+\t\leq T}\mE|U^m_{t+\t}-U^m_t|^p=0.
\label{equicont}
\ee

Finally, (\ref{umunbou}) and (\ref{equicont}) imply that $\{U^m\}_{m\in\mN}$ is tight in $C([0,T], \mathbb{R})$.

{\bf Step 2.} We show that the limit of $\{U^m\}_{m\in\mN}$ solves Eq.(\ref{ug}).

Let $\left\{U^{m_k}\right\}_{k\in\mN}$ be some weakly converging subsequence and $\tilde{U}$ be its limit point. Then we construct a weak solution for Eq.(\ref{ug}). By the Skorohod representation theorem, there exist a probability space $(\hat{\Omega},\hat{\sF}, \hat{\mP})$, two $C([0,T], \mathbb{R})$-valued random variables $\hat{U}, \hat{B}$ and a $C([0,T], \mathbb{R})$-valued random variable sequence $\{\hat{U}^{k}\}$ defined on it such that

(i) $(\hat{U}^{k},\hat{B}, \hat{U})\overset{d}{=}(U^{m_k},B,\tilde{U})$, for any $k\in\mN$;

(ii) $\hat{U}^{k}\overset{a.s.}{\rightarrow}\hat{U}$ as $k\rightarrow\infty$.

Set
$$
Y_t^k:=\int_0^t(t-s)^{\alpha-1} \hat{U}_s^{k} \dif s, \quad t\geq 0, \quad k\in\mN,
$$
and it is easy to check that $Y^k$ satisfies the following equation
$$
Y_t^k=\int_0^t(t-s)^{\alpha-1} g(s)\dif s+c_\alpha \int_0^t b^{m_k}\left(s,\hat{U}_s^{k}\right) \dif s+c_\alpha \int_0^t \sigma^{m_k}\left(s,\hat{U}_s^{k}\right) \dif \hat{B}_s,
$$
where $c_\alpha=\int_0^1(1-r)^{\alpha-1} r^{-\alpha}\dif r$. By passing to the limit, due to convergence of $\{\hat{U}^{k}\}_{k\in\mN}$ to $\hat{U}$  in $C([0,T], \mathbb{R})$ we get that 
\ce
Y_t^k=\int_0^t(t-s)^{\alpha-1} \hat{U}_s^{k} \dif s&\rightarrow& Y_t:=\int_0^t(t-s)^{\alpha-1} \hat{U}_s\dif s, \\
c_\alpha \int_0^t b^{m_k}\left(s,\hat{U}_s^{k}\right) \dif s&\rightarrow& c_\alpha \int_0^t b(s,\hat{U}_s)\dif s, \quad a.s. \hat{\mP}.
\de
Moreover, set
$$
M_t^k:=c_\alpha \int_0^t \sigma^{m_k}\left(s,\hat{U}_s^{k}\right) \dif \hat{B}_s, \quad t\in[0,T], \quad k \in\mN,
$$
and $M_{\cdot}^k$ is a sequence of square integrable martingales with quadratic variations given by
$$
\left\langle M_{\cdot}^k\right\rangle_t=c_\alpha^2 \int_0^t \sigma^{m_k}\left(s,\hat{U}_s^{k}\right)^2\dif s, \quad t\in[0,T], \quad k \in\mN.
$$
By uniform convergence of $\left\{\sigma^{m_k}\right\}_{k\in\mN}$ to $\sigma$ and convergence of $\{\hat{U}^{k}\}_{k\in\mN}$ to $\hat{U}$ in $C([0,T], \mathbb{R})$, we obtain that martingales $M_{\cdot}^k$  converge to the martingale 
$$
M_{\cdot}=c_\alpha \int_0^{\cdot} \sigma(s,\hat{U}_s)\dif \hat{B}_s, \quad t\in[0,T],
$$ 
and hence,
$$
Y_t=\int_0^t(t-s)^{\alpha-1}g(s)\dif s+c_\alpha \int_0^t b(s,\hat{U}_s)\dif s+c_\alpha \int_0^t \sigma(s,\hat{U}_s)\dif \hat{B}_s, \quad t\in[0,T].
$$
By reversing the transformation, that is, by recalling that
$$
\hat{U}_t=\frac{1}{c_\alpha} \frac{\dif }{\dif t} \int_0^t(t-s)^{-\alpha} Y_s \dif s, \quad t\in[0,T],
$$
it is easy to verify that $\hat{U}$ is a weak solution to Eq.(\ref{ug}).
\end{proof}

{\bf Proof of Theorem \ref{eu}.}

{\bf Existence.} Taking $g(t)=x_0$, by Proposition \ref{cont1} we obtain the existence for solutions of Eq.(\ref{Eq1}).

{\bf Uniqueness.} By the deduction below Remark \ref{spdeth} and Theorem \ref{exisuniq}, the uniqueness for solutions of Eq.(\ref{Eq1}) holds. The proof is complete.

\section{Proof of Theorem \ref{equimildweak} and \ref{exisuniq}}\label{thspdeeu}

In this section, we prove Theorem \ref{equimildweak} and \ref{exisuniq}. 

\subsection{Proof of Theorem \ref{equimildweak}}

In this subsection, we prove Theorem \ref{equimildweak}.

{\bf Weak solutions}$\Longrightarrow${\bf Mild solutions.} If ${\bf X}$ is a weak solution to Eq.(\ref{spde}), it holds that for any $\phi_{\cdot}(\cdot)\in C([0,T], \cD(\Delta_\t))$ and $s\rightarrow\frac{\dif \phi_s}{\dif s}\in C([0,T], C_{rap})$,
\ce
\<{\bf X}_t, \phi_t\>&=&\<{\bf X}_0, \phi_0\>+\int_0^t\<{\bf X}_s,\frac{\dif \phi_s}{\dif s}+\Delta_\t\phi_s\>\dif s+\int_0^t\<\frac{1}{c_\t}b(s,{\bf X}_s(0))\d_0,\phi_s\>\dif s\\
&&+\int_0^t\<\frac{1}{c_\t}\sigma(s,{\bf X}_s(0))\d_0\dif B_s,\phi_s\>.
\de
Then we take $\phi_s=S^{\t}_{t-s}\phi$ for $\phi\in\cD(\Delta_\t)$ and obtain that
\ce
\<{\bf X}_t, \phi\>&=&\<{\bf X}_0, S^{\t}_{t}\phi\>+\int_0^t\<{\bf X}_s,\frac{\dif S^{\t}_{t-s}\phi}{\dif s}+\Delta_\t S^{\t}_{t-s}\phi\>\dif s+\int_0^t\<\frac{1}{c_\t}b(s,{\bf X}_s(0))\d_0,S^{\t}_{t-s}\phi\>\dif s\\
&&+\int_0^t\<\frac{1}{c_\t}\sigma(s,{\bf X}_s(0))\d_0\dif B_s,S^{\t}_{t-s}\phi\>.
\de
Note that $\Delta_\t$ is the generator of $\{S^{\t}_{t}\}$ and
$$
\frac{\dif S^{\t}_{t-s}\phi}{\dif s}+\Delta_\t S^{\t}_{t-s}\phi=-\frac{\dif S^{\t}_{t-s}\phi}{\dif (t-s)}+\Delta_\t S^{\t}_{t-s}\phi=0.
$$
Thus, some computation implies that
\ce
\<{\bf X}_t, \phi\>&=&\<S^{\t}_{t}{\bf X}_0, \phi\>+\int_0^t\<S^{\t}_{t-s}\frac{1}{c_\t}b(s,{\bf X}_s(0))\d_0,\phi\>\dif s+\int_0^t\<S^{\t}_{t-s}\frac{1}{c_\t}\sigma(s,{\bf X}_s(0))\d_0\dif B_s,\phi\>,
\de
which yields that
\ce
{\bf X}_t=S_t^\t {\bf X}_0+\int_0^tS_{t-s}^\t \frac{1}{c_\t}b(s,{\bf X}_s(0))\d_0\dif s+\int_0^tS_{t-s}^\t\frac{1}{c_\t}\sigma(s,{\bf X}_s(0))\d_0\dif B_s.
\de

{\bf Mild solutions}$\Longrightarrow${\bf Weak solutions.} If ${\bf X}$ is a mild solution to Eq.(\ref{spde}), it holds that
\ce
{\bf X}_t=S_t^\t {\bf X}_0+\int_0^tS_{t-s}^\t \frac{1}{c_\t}b(s,{\bf X}_s(0))\d_0\dif s+\int_0^tS_{t-s}^\t\frac{1}{c_\t}\sigma(s,{\bf X}_s(0))\d_0\dif B_s.
\de
Then for $\phi\in\cD(\Delta_\t)$, we have that
\be
\int_0^t\<{\bf X}_s,\Delta_\t\phi\>\dif s&=&\int_0^t\<S_s^\t {\bf X}_0,\Delta_\t\phi\>\dif s+\int_0^t\<\int_0^sS_{s-r}^\t \frac{1}{c_\t}b(r,{\bf X}_r(0))\d_0\dif r,\Delta_\t\phi\>\dif s\no\\
&&+\int_0^t\<\int_0^sS_{s-r}^\t \frac{1}{c_\t}\s(r,{\bf X}_r(0))\d_0\dif B_r,\Delta_\t\phi\>\dif s\no\\
&=:&I_1+I_2+I_3.
\label{i123}
\ee

First of all, for $I_1$, it holds that
\be
I_1=\int_0^t\<\Delta_\t S_s^\t {\bf X}_0,\phi\>\dif s=\int_0^t\<\frac{\dif S_s^\t {\bf X}_0}{\dif s},\phi\>\dif s=\<S_t^\t {\bf X}_0,\phi\>-\<{\bf X}_0,\phi\>.
\label{i1}
\ee
Then we deal with $I_2$ and obtain that
\be
I_2&=&\int_0^t\dif r\int_r^t \<S_{s-r}^\t\frac{1}{c_\t}b(r,{\bf X}_r(0))\d_0,\Delta_\t\phi\>\dif s\no\\
&=&\int_0^t\dif r\int_r^t \<\Delta_\t S_{s-r}^\t\frac{1}{c_\t}b(r,{\bf X}_r(0))\d_0,\phi\>\dif s\no\\
&=&\int_0^t\dif r\int_r^t \<\frac{\dif S_{s-r}^\t}{\dif (s-r)}\frac{1}{c_\t}b(r,{\bf X}_r(0))\d_0,\phi\>\dif s\no\\
&=&\int_0^t \<S_{t-r}^\t\frac{1}{c_\t}b(r,{\bf X}_r(0))\d_0,\phi\>\dif r-\int_0^t \<\frac{1}{c_\t}b(r,{\bf X}_r(0))\d_0,\phi\>\dif r.
\label{i2}
\ee
By the same deduction to that for $I_2$, it follows that
\be
I_3=\int_0^t \<S_{t-r}^\t\frac{1}{c_\t}\s(r,{\bf X}_r(0))\d_0\dif B_r,\phi\>-\int_0^t \<\frac{1}{c_\t}b(r,{\bf X}_r(0))\d_0\dif B_r,\phi\>.
\label{i3}
\ee

Finally, by combining (\ref{i1}) (\ref{i2}) (\ref{i3}) with (\ref{i123}),  we conclude that 
\ce
\<{\bf X}_t,\phi\>&=&\<{\bf X}_0,\phi\>+\int_0^t\<{\bf X}_s,\Delta_\t\phi\>\dif s+\int_0^t \<\frac{1}{c_\t}b(r,{\bf X}_r(0))\d_0,\phi\>\dif r\\
&&+\int_0^t \<\frac{1}{c_\t}b(r,{\bf X}_r(0))\d_0\dif B_r,\phi\>.
\de
The proof is complete.

\subsection{Proof of Theorem \ref{exisuniq}}

In this subsection, we prove Theorem \ref{exisuniq}. And we show the existence and uniqueness of mild solutions to Eq.(\ref{spde}), respectively. 

{\bf Existence.} First of all, we consider the following equation:
\be
\bar{U}_t=\int_{\mR}p_{t}^\t(0-y){\bf X}_0(y)\dif y+\int_0^t(t-s)^{-\a}b(s,\bar{U}_s)\dif s+\int_0^t(t-s)^{-\a}\sigma(s,\bar{U}_s)\dif B_s.
\label{ugpt}
\ee
And we claim that $\int_{\mR}p_{t}^\t(0-y){\bf X}_0(y)\dif y$ is continuous in $t$. In fact, for any $T>0$ and $t_1, t_2\in[0,T]$, it holds that
\ce
&&\left|\int_{\mR}p_{t_2}^\t(0-y){\bf X}_0(y)\dif y-\int_{\mR}p_{t_1}^\t(0-y){\bf X}_0(y)\dif y\right|\\
&=&c_\t\left|\int_{\mR}\frac{e^{-\frac{|y|^{2+\t}}{2t_2}}}{t_2^\a}{\bf X}_0(y)\dif y-\int_{\mR}\frac{e^{-\frac{|y|^{2+\t}}{2t_1}}}{t_1^\a}{\bf X}_0(y)\dif y\right|\\
&=&c_\t\left|\int_{\mR}e^{-\frac{|z|^{2+\t}}{2}}{\bf X}_0(t_2^\a z)\dif z-\int_{\mR}e^{-\frac{|z|^{2+\t}}{2}}{\bf X}_0(t_1^\a z)\dif z\right|\\
&\leq&c_\t\int_{\mR}e^{-\frac{|z|^{2+\t}}{2}}\left|{\bf X}_0(t_2^\a z)-{\bf X}_0(t_1^\a z)\right|\dif z.
\de
Since ${\bf X}_0\in C_{tem}$, we obtain that
$$
\lim\limits_{t_2\rightarrow t_1}\left|{\bf X}_0(t_2^\a z)-{\bf X}_0(t_1^\a z)\right|=0,
$$
and
\ce
&&e^{-\frac{|z|^{2+\t}}{2}}\left|{\bf X}_0(t_2^\a z)-{\bf X}_0(t_1^\a z)\right|\\
&\leq&e^{-\frac{|z|^{2+\t}}{2}+\lambda T^\a |z|}\left(\sup\limits_{t_2^\a z\in\mR}|{\bf X}_0(t_2^\a z)|e^{-\lambda |t_2^\a z|}+\sup\limits_{t_1^\a z\in\mR}|{\bf X}_0(t_1^\a z)|e^{-\lambda |t_1^\a z|}\right)\\
&\leq&2e^{-\frac{|z|^{2+\t}}{2}+\lambda T^\a |z|}\|{\bf X}_0\|_{\lambda,\infty}.
\de
Noting that for suitable $\lambda>0$, $\int_{\mR}e^{-\frac{|z|^{2+\t}}{2}+\lambda T^\a |z|}\dif z<\infty$. Therefore, the dominated convergence theorem implies that $\int_{\mR}p_{t}^\t(0-y){\bf X}_0(y)\dif y$ is continuous in $t$. Thus, by Proposition \ref{cont1}, we know that Eq.(\ref{ugpt}) has a solution denoted as $\bar{U}$.

Set
\ce
{\bf X}_t(x)&:=&\int_{\mR}p_{t}^\t(x-y){\bf X}_0(y)\dif y+\int_0^t\int_{\mR}p_{t-s}^\t(x-y)\frac{b(s,\bar{U}_s)}{c_\t}\d_0(y)\dif y\dif s\\
&&+\int_0^t\int_{\mR}p_{t-s}^\t(x-y)\frac{\sigma(s,\bar{U}_s)}{c_\t}\d_0(y)\dif y\dif B_s,
\de
and it is easy to check that ${\bf X}$ is a mild solution to Eq.(\ref{spde}) with ${\bf X}_t(0)=\bar{U}_t, t\geq 0$, i.e.
\ce
{\bf X}_t(x)&=&(S_t^\t {\bf X}_0)(x)+\int_0^t\left(S_{t-s}^\t\frac{b(s,{\bf X}_s(0))}{c_\t}\d_0\right)(x)\dif s\\
&&+\int_0^t\left(S_{t-s}^\t\frac{\sigma(s,{\bf X}_s(0))}{c_\t}\d_0\right)(x)\dif B_s,
\de
 and ${\bf X}\in C(\mR_+, C_{tem})$. Moreover, by the simple calculation and (\ref{ubou}), it holds that for any $T, \lambda>0$
\be
\sup\limits_{0\leq t\leq T}\sup\limits_{x\in\mR}\mE\left(|{\bf X}_t(x)|^pe^{-\lambda |x|}\right)<\infty, \quad p>\frac{2}{1-2\a}.
\label{tlamp}
\ee

{\bf Uniqueness.} The method below is taken from \cite{mps}. 

{\bf Step 1.} Let ${\bf X}^1, {\bf X}^2$ be two weak solutions of Eq.(\ref{spde}) with ${\bf X}^1_0={\bf X}^2_0={\bf X}_0$. Set ${\bf Z}_t:={\bf X}^1_t-{\bf X}^2_t$, and for $\Phi_x^m(y)=p^\t_{m^{-1/\a}}(x-y), m\in\mN$, by (\ref{weaksoluspde}) it holds that
\ce
\<{\bf Z}_t,\Phi_x^m\>&=&\int_0^t\left<{\bf Z}_r,\Delta_\t\Phi_x^m\right>\dif r+\int_0^t\frac{1}{c_\t}\(b(r,{\bf X}^1_r(0))-b(r,{\bf X}^2_r(0))\)\Phi_x^m(0)\dif r\\
&&+\int_0^t\frac{1}{c_\t}\(\sigma(r,{\bf X}^1_r(0))-\sigma(r,{\bf X}^2_r(0))\)\Phi_x^m(0)\dif B_r.
\de

Next, we define a sequence of functions $\phi_n$ as follows. First of all, let $\{a_n\}$ be a strictly decreasing sequence such that $a_0=1, \lim\limits_{n\rightarrow\infty}a_n=0$ and
$$
\int_{a_n}^{a_{n-1}}\frac{1}{x}\dif x=n, \quad n\geq 1, \quad n\in\mN.
$$
Then let $\{\psi_n\}$ be continuous functions such that ${\rm supp}(\psi_n)\subset(a_n, a_{n-1})$ and 
\be
0\leq \psi_n(x)\leq \frac{2}{nx}, ~\mbox{for}~\forall x>0, ~\mbox{and}~ \int_{a_n}^{a_{n-1}}\psi_n(x)\dif x=1.
\label{psibou}
\ee
Finally, set 
$$
\phi_n(x):=\int_0^{|x|}\int_0^y\psi_n(z)\dif z\dif y,
$$
and we know that $\phi_n(x)\uparrow |x|$ and 
\ce
&&\phi_n^{\prime}(x)={\rm sgn}(x)\int_0^{|x|}\psi_n(y)\dif y,\\
&&\phi_n^{\prime\prime}(x)=\psi_n(|x|).
\de
Moreover, $|\phi_n^{\prime}(x)|\leq 1$ and as $n\rightarrow\infty$, $\int_{\mR} \phi_n^{\prime\prime}(x)h(x)\dif x\rightarrow h(0)$ for any function $h$ which is continuous  at $0$.

In the following, applying the It\^o formula to $\phi_n(\<{\bf Z}_t,\Phi_x^m\>)$, we obtain that
\be
\phi_n(\<{\bf Z}_t,\Phi_x^m\>)&=&\int_0^t\phi'_n(\<{\bf Z}_r,\Phi_x^m\>)\left<{\bf Z}_r,\Delta_\t\Phi_x^m\right>\dif r\no\\
&&+\int_0^t\phi'_n(\<{\bf Z}_r,\Phi_x^m\>)\frac{1}{c_\t}\(b(r,{\bf X}^1_r(0))-b(r,{\bf X}^2_r(0))\)\Phi_x^m(0)\dif r\no\\
&&+\int_0^t\phi'_n(\<{\bf Z}_r,\Phi_x^m\>)\frac{1}{c_\t}\(\sigma(r,{\bf X}^1_r(0))-\sigma(r,{\bf X}^2_r(0))\)\Phi_x^m(0)\dif B_r\no\\
&&+\frac{1}{2}\int_0^t\phi''_n(\<{\bf Z}_r,\Phi_x^m\>)\frac{1}{c^2_\t}\(\sigma(r,{\bf X}^1_r(0))-\sigma(r,{\bf X}^2_r(0))\)^2\Phi_x^m(0)^2\dif r.
\label{phiPhi}
\ee
Take a nonnegative test function $\Psi\in C([0,\infty), \cD(\Delta_\t))$ such that 
\be
\Psi_s(0)>0, \forall s\geq 0 ~\mbox{and}~ \sup\limits_{0\leq s\leq t}\left|\int_{0^+}|x|^{-\t}\left(\frac{\partial \Psi_s(x)}{\partial x}\right)^2(\Psi_s(x))^{-1}\dif x\right|<\infty, \quad \forall t>0,
\label{Psipro}
\ee
and $s\mapsto\frac{\partial \Psi_s(\cdot)}{\partial s}\in C(\mR_+, C_{rap})$. Also assume that $\Gamma(t):=\{x: \exists s\leq t, \Psi_s(x)>0\}\subset B(0, J(t))$ for some $J(t)>0$. So, applying the It\^o formula to $\phi_n(\<{\bf Z}_t,\Phi_x^m\>)\Psi_t(x)$, integrating two sides with respect to $x$ and taking the expectation on two sides, we have that
\be
&&\mE\<\phi_n(\<{\bf Z}_t,\Phi_{\cdot}^m\>),\Psi_t\>\no\\
&=&\mE\int_0^t\<\phi'_n(\<{\bf Z}_r,\Phi_{\cdot}^m\>)\left<{\bf Z}_r,\Delta_\t\Phi_{\cdot}^m\right>,\Psi_r\>\dif r\no\\
&&+\mE\int_0^t\<\phi'_n(\<{\bf Z}_r,\Phi_{\cdot}^m\>)\Phi_{\cdot}^m(0),\Psi_r\>\frac{1}{c_\t}\(b(r,{\bf X}^1_r(0))-b(r,{\bf X}^2_r(0))\)\dif r\no\\
&&+\mE\int_0^t\<\phi'_n(\<{\bf Z}_r,\Phi_{\cdot}^m\>)\Phi_{\cdot}^m(0),\Psi_r\>\frac{1}{c_\t}\(\sigma(r,{\bf X}^1_r(0))-\sigma(r,{\bf X}^2_r(0))\)\dif B_r\no\\
&&+\frac{1}{2}\mE\int_0^t\<\phi''_n(\<{\bf Z}_r,\Phi_{\cdot}^m\>)\Phi_{\cdot}^m(0)^2,\Psi_r\>\frac{1}{c^2_\t}\(\sigma(r,{\bf X}^1_r(0))-\sigma(r,{\bf X}^2_r(0))\)^2\dif r\no\\
&&+\mE\int_0^t\<\phi_n(\<{\bf Z}_r,\Phi_{\cdot}^m\>),\frac{\partial \Psi_r(\cdot)}{\partial r}\>\dif r.
\label{i1i2i3i4i5}
\ee
We claim that
\be
\lim\limits_{n,m\rightarrow\infty}\mE\<\phi_n(\<{\bf Z}_t,\Phi_{\cdot}^m\>),\Psi_t\>\leq\mE\int_0^t\int_{\mR}|{\bf Z}_r(x)|\left(\Delta_\t\Psi_r(x)+\frac{\partial \Psi_r(x)}{\partial r}\right)\dif x\dif r. 
\label{clai}
\ee
Based on the claim and the Fatou lemma, it holds that
\be
\int_{\mR}\mE|{\bf Z}_t(x)| \Psi_t(x)\dif x\leq \int_0^t\int_{\mR}\mE|{\bf Z}_r(x)|\left|\Delta_\t\Psi_r(x)+\frac{\partial \Psi_r(x)}{\partial r}\right|\dif x\dif r.
\label{zrli}
\ee

Next, let $\{g_N, N\in\mN\}$ be a sequence of functions in $C_c^\infty(\mR)$ such that $g_N: \mR\mapsto [0,1]$, $g_N(x)=1$ for $|x|\leq N$, $g_N(x)=0$ for $|x|>N+1$ and 
$$
\sup\limits_{N}\(\||x|^{-\t}g'_N\|_{\infty}+\|\Delta_\t g_N\|_{\infty}\)<\infty.
$$
Again set for $(r,x)\in[0,t]\times\mR$
$$
\Psi^{(N)}_r(x):=(S^\t_{t-r}h)(x)g_N(x), \quad h\in C_c^\infty(\mR),
$$
and one can verify that $\Psi^{(N)}_r\in C_c^\infty(\mR)$ and for $\lambda>0$, there exists a constant $C>0$ such that for all $N$
\ce
\left|\Delta_\t \Psi^{(N)}_r(x)+\frac{\partial \Psi^{(N)}_r(x)}{\partial r}\right|&=&\left|4\a^2|x|^{-\t}\frac{\partial (S^\t_{t-r}h)(x)}{\partial x}\frac{\partial g_N(x)}{\partial x}+(S^\t_{t-r}h)(x)\Delta_\t g_N(x)\right|\\
&\leq& Ce^{-\lambda|x|}I_{|x|>N}.
\de
Replacing $\Psi_t$ by $\Psi^{(N)}_t$ in (\ref{zrli}), by the above inequality, we conclude that
$$
\int_{\mR}\mE|{\bf Z}_t(x)| h(x)g_N(x)\dif x\leq C\int_0^t\int_{\mR}\mE|{\bf Z}_r(x)|e^{-\lambda|x|}I_{|x|>N}\dif x\dif r.
$$
Note that by (\ref{tlamp}), the right side of the above inequality tends to zero as $N\rightarrow\infty$. Therefore, for $t\geq 0, x\in\mR$
\ce
\mE|{\bf Z}_t(x)|=0,
\de
and furthermore
$$
{\bf X}^1_t(x)={\bf X}^2_t(x).
$$
Finally, the continuity of ${\bf X}^1_t(x), {\bf X}^2_t(x)$ in $t, x$ gives the required result.

{\bf Step 2.} We prove the claim (\ref{clai}). To do this, by (\ref{i1i2i3i4i5}) we divide $\mE\<\phi_n(\<{\bf Z}_t,\Phi_{\cdot}^m\>),\Psi_t\>$ into $I_1, I_2, I_3, I_4, I_5$, where
\ce
&&I_1:=\mE\int_0^t\<\phi'_n(\<{\bf Z}_r,\Phi_{\cdot}^m\>)\left<{\bf Z}_r,\Delta_\t\Phi_{\cdot}^m\right>,\Psi_r\>\dif r,\\
&&I_2:=\mE\int_0^t\<\phi'_n(\<{\bf Z}_r,\Phi_{\cdot}^m\>)\Phi_{\cdot}^m(0),\Psi_r\>\frac{1}{c_\t}\(b(r,{\bf X}^1_r(0))-b(r,{\bf X}^2_r(0))\)\dif r,\\
&&I_3:=\mE\int_0^t\<\phi'_n(\<{\bf Z}_r,\Phi_{\cdot}^m\>)\Phi_{\cdot}^m(0),\Psi_r\>\frac{1}{c_\t}\(\sigma(r,{\bf X}^1_r(0))-\sigma(r,{\bf X}^2_r(0))\)\dif B_r,\\
&&I_4:=\frac{1}{2}\mE\int_0^t\<\phi''_n(\<{\bf Z}_r,\Phi_{\cdot}^m\>)\Phi_{\cdot}^m(0)^2,\Psi_r\>\frac{1}{c^2_\t}\(\sigma(r,{\bf X}^1_r(0))-\sigma(r,{\bf X}^2_r(0))\)^2\dif r,\\
&&I_5:=\mE\int_0^t\<\phi_n(\<{\bf Z}_r,\Phi_{\cdot}^m\>),\frac{\partial \Psi_r(\cdot)}{\partial r}\>\dif r.
\de

{\bf Estimation for $I_1$.} Note that
$$
\int_{\mR}{\bf Z}_r(y)\Delta_{y,\t}\Phi_x^m(y)\dif y=\int_{\mR}{\bf Z}_r(y)\Delta_{x,\t}\Phi_x^m(y)\dif y=\Delta_{x,\t}\int_{\mR}{\bf Z}_r(y)\Phi_x^m(y)\dif y,
$$
where the fact $\Delta_{y,\t}\Phi_x^m(y)=\Delta_{x,\t}\Phi_x^m(y)$ is used. Thus, by integration by parts, it holds that
\ce
&&\int_0^t\<\phi'_n(\<{\bf Z}_r,\Phi_{\cdot}^m\>)\left<{\bf Z}_r,\Delta_\t\Phi_{\cdot}^m\right>,\Psi_r\>\dif r\\
&=&\int_0^t\int_{\mR}\phi'_n(\<{\bf Z}_r,\Phi_x^m\>)\Delta_{x,\t}(\<{\bf Z}_r,\Phi_x^m\>)\Psi_r(x)\dif x\dif r\\
&=&-2\a^2\int_0^t\int_{\mR}\frac{\partial}{\partial x}\(\phi'_n(\<{\bf Z}_r,\Phi_x^m\>)\)|x|^{-\t}\frac{\partial}{\partial x}\(\<{\bf Z}_r,\Phi_x^m\>\)\Psi_r(x)\dif x\dif r\\
&&-2\a^2\int_0^t\int_{\mR}\phi'_n(\<{\bf Z}_r,\Phi_x^m\>)|x|^{-\t}\frac{\partial}{\partial x}\(\<{\bf Z}_r,\Phi_x^m\>\)\frac{\partial}{\partial x}\Psi_r(x)\dif x\dif r\\
&=&-2\a^2\int_0^t\int_{\mR}\phi''_n(\<{\bf Z}_r,\Phi_x^m\>)|x|^{-\t}\(\frac{\partial}{\partial x}\(\<{\bf Z}_r,\Phi_x^m\>\)\)^2\Psi_r(x)\dif x\dif r\\
&&+2\a^2\int_0^t\int_{\mR}\phi''_n(\<{\bf Z}_r,\Phi_x^m\>)|x|^{-\t}\frac{\partial}{\partial x}\(\<{\bf Z}_r,\Phi_x^m\>\)\<{\bf Z}_r,\Phi_x^m\>\frac{\partial}{\partial x}\Psi_r(x)\dif x\dif r\\
&&+\int_0^t\int_{\mR}\phi'_n(\<{\bf Z}_r,\Phi_x^m\>)\<{\bf Z}_r,\Phi_x^m\>\Delta_\t\Psi_r(x)\dif x\dif r\\
&=&2\a^2\int_0^t\int_{\mR}\phi''_n(\<{\bf Z}_r,\Phi_x^m\>)|x|^{-\t}\[-\(\frac{\partial}{\partial x}\(\<{\bf Z}_r,\Phi_x^m\>\)\)^2\Psi_r(x)\\
&&\qquad +\frac{\partial}{\partial x}\(\<{\bf Z}_r,\Phi_x^m\>\)\<{\bf Z}_r,\Phi_x^m\>\frac{\partial}{\partial x}\Psi_r(x)\]\dif x\dif r\\
&&+\int_0^t\int_{\mR}\phi'_n(\<{\bf Z}_r,\Phi_x^m\>)\<{\bf Z}_r,\Phi_x^m\>\Delta_\t\Psi_r(x)\dif x\dif r\\
&=&J_{11}+J_{12}.
\de
For any $r\in[0,t]$, set
\ce
A^r&:=&\left\{x: \(\frac{\partial}{\partial x}\(\<{\bf Z}_r,\Phi_x^m\>\)\)^2\Psi_r(x)\leq \frac{\partial}{\partial x}\(\<{\bf Z}_r,\Phi_x^m\>\)\<{\bf Z}_r,\Phi_x^m\>\frac{\partial}{\partial x}\Psi_r(x)\right\}\cap\{x: \Psi_r(x)>0\}\\
&=&A^{+,r}\cup A^{-,r}\cup A^{0,r},
\de
where
\ce
&&A^{+,r}:=A^r\cap\left\{x: \frac{\partial}{\partial x}\(\<{\bf Z}_r,\Phi_x^m\>\)>0\right\},\\
&&A^{-,r}:=A^r\cap\left\{x: \frac{\partial}{\partial x}\(\<{\bf Z}_r,\Phi_x^m\>\)<0\right\},\\
&&A^{0,r}:=A^r\cap\left\{x: \frac{\partial}{\partial x}\(\<{\bf Z}_r,\Phi_x^m\>\)=0\right\}.
\de

Note that on $A^{+,r}$, it holds that
$$
0<\frac{\partial}{\partial x}\(\<{\bf Z}_r,\Phi_x^m\>\)\Psi_r(x)\leq\<{\bf Z}_r,\Phi_x^m\>\frac{\partial}{\partial x}\Psi_r(x),
$$
which implies that
\ce
&&\int_0^t\int_{A^{+,r}}\psi_n(|\<{\bf Z}_r,\Phi_x^m\>|)|x|^{-\t}\frac{\partial}{\partial x}\(\<{\bf Z}_r,\Phi_x^m\>\)\<{\bf Z}_r,\Phi_x^m\>\frac{\partial}{\partial x}\Psi_r(x)\dif x\dif r\\
&\leq&\int_0^t\int_{A^{+,r}}\psi_n(|\<{\bf Z}_r,\Phi_x^m\>|)|x|^{-\t}\<{\bf Z}_r,\Phi_x^m\>^2\frac{(\frac{\partial}{\partial x}\Psi_r(x))^2}{\Psi_r(x)}\dif x\dif r\\
&\leq&\int_0^t\int_{A^{+,r}}\frac{2}{n}I_{\{a_{n}\leq |\<{\bf Z}_r,\Phi_x^m\>|\leq a_{n-1}\}}|x|^{-\t}|\<{\bf Z}_r,\Phi_x^m\>|\frac{(\frac{\partial}{\partial x}\Psi_r(x))^2}{\Psi_r(x)}\dif x\dif r\\
&\leq&\frac{2a_{n-1}}{n}\int_0^t\int_{\mR}I_{\{\Psi_r(x)>0\}}|x|^{-\t}\frac{(\frac{\partial}{\partial x}\Psi_r(x))^2}{\Psi_r(x)}\dif x\dif r\\
&\leq&\frac{2a_{n-1}}{n}\int_0^t\(\int_{B(0,\e)}|x|^{-\t}\frac{(\frac{\partial}{\partial x}\Psi_r(x))^2}{\Psi_r(x)}\dif x+\int_{\Gamma(t)\setminus B(0,\e)}|x|^{-\t}\frac{(\frac{\partial}{\partial x}\Psi_r(x))^2}{\Psi_r(x)}\dif x\)\dif r\\
&\overset{(\ref{Psipro})}{\leq}&\frac{2a_{n-1}}{n}\int_0^t\(C+2\left\|\frac{\partial^2}{\partial x^2}\Psi_r\right\|_{\infty}\int_{\Gamma(t)\setminus B(0,\e)}|x|^{-\t}\dif x\)\dif r\\
&\leq&\frac{2a_{n-1}}{n}C_{\Psi,t},
\de
where Lemma \ref{psiseco} below is used in the last second inequality and the constant $\e>0$ is small enough such that 
$$
B(0,\e)\subset \Gamma(t), ~\mbox{and}~ \inf\limits_{r\leq t, x\in B(0,\e)}\Psi_r(x)>0.
$$

On $A^{-,r}$, we know that
\ce
&&0>\frac{\partial}{\partial x}\(\<{\bf Z}_r,\Phi_x^m\>\)\Psi_r(x)\geq\<{\bf Z}_r,\Phi_x^m\>\frac{\partial}{\partial x}\Psi_r(x),\\
&&\frac{\partial}{\partial x}\(\<{\bf Z}_r,\Phi_x^m\>\)\<{\bf Z}_r,\Phi_x^m\>\frac{\partial}{\partial x}\Psi_r(x)\leq\<{\bf Z}_r,\Phi_x^m\>^2\frac{(\frac{\partial}{\partial x}\Psi_r(x))^2}{\Psi_r(x)}.
\de
Therefore, by similar deduction to the above, one can get that
\ce
\int_0^t\int_{A^{-,r}}\psi_n(|\<{\bf Z}_r,\Phi_x^m\>|)|x|^{-\t}\frac{\partial}{\partial x}\(\<{\bf Z}_r,\Phi_x^m\>\)\<{\bf Z}_r,\Phi_x^m\>\frac{\partial}{\partial x}\Psi_r(x)\dif x\dif r\leq\frac{2a_{n-1}}{n}C_{\Psi,t}.
\de

Finally, it is easy to see that 
$$
\int_0^t\int_{A^{0,r}}\psi_n(|\<{\bf Z}_r,\Phi_x^m\>|)|x|^{-\t}\frac{\partial}{\partial x}\(\<{\bf Z}_r,\Phi_x^m\>\)\<{\bf Z}_r,\Phi_x^m\>\frac{\partial}{\partial x}\Psi_r(x)\dif x\dif r=0.
$$

By the above deduction, we obtain that
\ce
\mE J_{11}\leq 4\a^2C_{\Psi,t}\frac{a_{n-1}}{n}.
\de

Now we observe $J_{12}$. Note that $\phi'_n(x)x\uparrow |x|$ uniformly in $x$ as $n\rightarrow\infty$, and $\<{\bf Z}_r,\Phi_x^m\>$ tends to ${\bf Z}_r(x)$ a.s. as $m\rightarrow\infty$ for all $r,x$. Hence, $\phi'_n(\<{\bf Z}_r,\Phi_x^m\>)\<{\bf Z}_r,\Phi_x^m\>\rightarrow|{\bf Z}_r(x)|$ a.s. as $n,m\rightarrow\infty$. Besides, the Jensen inequality and (\ref{tlamp}) imply that $\<|{\bf Z}_r|,\Phi_x^m\>$ is $L^p$ bounded on 
$([0,t]\times B(0,J(t))\times \Omega, \dif r\times\dif x\times \mP)$ uniformly in $m$. This gives uniform integrability of $\<|{\bf Z}_r|,\Phi_x^m\>$ in $m$ on $[0,t]\times B(0,J(t))\times \Omega$. Moreover, we know that
$$
|\phi'_n(\<{\bf Z}_r,\Phi_x^m\>)\<{\bf Z}_r,\Phi_x^m\>|\leq |\<{\bf Z}_r,\Phi_x^m\>|\leq \<|{\bf Z}_r|,\Phi_x^m\>.
$$
So, $\{\phi'_n(\<{\bf Z}_r,\Phi_x^m\>)\<{\bf Z}_r,\Phi_x^m\>\}$ is uniformly integrable. Since $\Psi_r=0$ outside $B(0,J(t))$, this yields that
\ce
\lim\limits_{n,m\rightarrow\infty}\mE J_{12}=\mE\int_0^t\int_{\mR}|{\bf Z}_r(x)|\Delta_\t\Psi_r(x)\dif x\dif r.
\de

Collecting the pieces, we obtain that
\be
\lim\limits_{n,m\rightarrow\infty}I_{1}=\mE\int_0^t\int_{\mR}|{\bf Z}_r(x)|\Delta_\t\Psi_r(x)\dif x\dif r.
\label{i1}
\ee

{\bf Estimation for $I_2$.} By the Fatou lemma and (${\bf H}^1_{b}$), it holds that
\be
\limsup\limits_{m\rightarrow\infty}I_2&\leq&\mE\int_0^t\phi'_n({\bf Z}_r(0))\Psi_r(0)\frac{1}{c_\t}\(b(r,{\bf X}^1_r(0))-b(r,{\bf X}^2_r(0))\)\dif r\no\\
&\leq&0.
\label{i2}
\ee

{\bf Estimation for $I_3$.} Note that the quadratic variation of 
$$
\int_0^t\<\phi'_n(\<{\bf Z}_r,\Phi_{\cdot}^m\>)\Phi_{\cdot}^m(0),\Psi_r\>\frac{1}{c_\t}\(\sigma(r,{\bf X}^1_r(0))-\sigma(r,{\bf X}^2_r(0))\)\dif B_r
$$ 
satisfies
\ce
&&\int_0^t\<\phi'_n(\<{\bf Z}_r,\Phi_{\cdot}^m\>)\Phi_{\cdot}^m(0),\Psi_r\>^2\frac{1}{c^2_\t}|\sigma(r,{\bf X}^1_r(0))-\sigma(r,{\bf X}^2_r(0))|^2\dif r\\
&\leq&4L^2_\s\int_0^t\<\phi'_n(\<{\bf Z}_r,\Phi_{\cdot}^m\>)\Phi_{\cdot}^m(0),\Psi_r\>^2\frac{1}{c^2_\t}(2+|{\bf X}^1_r(0)|^2+|{\bf X}^2_r(0)|^2)\dif r\\
&\leq&4L^2_\s\|\Psi\|^2_\infty\frac{1}{c^2_\t}\int_0^t(2+|{\bf X}^1_r(0)|^2+|{\bf X}^2_r(0)|^2)\dif r.
\de
Thus, for any $t\in [0,T]$, we have that
\be
I_3=0.
\label{i3}
\ee

{\bf Estimation for $I_4$.} Note that $\<{\bf Z}_r,\Phi_x^m\>$ tends to ${\bf Z}_r(x)$ a.s. as $m\rightarrow\infty$ for all $r$ and $x$, and that $\Phi_{\cdot}^m(0)$ 
converges weakly to $\d_0$ as $m\rightarrow\infty$. So, we have that
\ce
\lim\limits_{m\rightarrow\infty}\<\phi''_n(\<{\bf Z}_r,\Phi_{\cdot}^m\>)\Phi_{\cdot}^m(0)^2,\Psi_r\>=\psi_n(|{\bf Z}_r(0)|)\Psi_r(0),
\de
which together with the Fatou lemma, (${\bf H}^1_{\sigma}$) and (\ref{psibou}) implies that
\ce
\limsup\limits_{m\rightarrow\infty}I_4&\leq&\frac{1}{2c^2_\t}\mE\int_0^t\psi_n(|{\bf Z}_r(0)|)\Psi_r(0)\(\sigma(r,{\bf X}^1_r(0))-\sigma(r,{\bf X}^2_r(0))\)^2\dif r\\
&\leq&\frac{L^2_\s}{nc^2_\t}\mE\int_0^t\Psi_r(0)|{\bf Z}_r(0)|^{2\g_2-1}\dif r\\
&\leq&\frac{C}{n},
\de
where (\ref{tlamp}) is used in the last inequality. Thus, we obtain that
\be
\lim\limits_{n,m\rightarrow\infty}I_{4}=0.
\label{i4}
\ee

{\bf Estimation for $I_5$.} Note that
$$
\phi_n(\<{\bf Z}_r,\Phi_x^m\>)\leq |\<{\bf Z}_r,\Phi_x^m\>|\leq \<|{\bf Z}_r|,\Phi_x^m\>,
$$
which together with the uniform integrability of $\<|{\bf Z}_r|,\Phi_x^m\>$ in $m$ on $[0,t]\times B(0,J(t))\times \Omega$ implies that $\{\phi_n(\<{\bf Z}_r,\Phi_{\cdot}^m\>): n,m\}$ is uniformly integrable on $[0,t]\times B(0,J(t))\times \Omega$. Moreover, we know that
$\phi_n(\<{\bf Z}_r,\Phi_x^m\>)\rightarrow|{\bf Z}_r(x)|$ as $n,m\rightarrow\infty$ a.s. for all $x$ and all $r\leq t$. Thus, it follows that
\be
\lim\limits_{n,m\rightarrow\infty}I_{5}=\mE\int_0^t\int_{\mR}|{\bf Z}_r(x)|\frac{\partial \Psi_r(x)}{\partial r}\dif x\dif r,
\label{i5}
\ee 
where we use the fact that
$$
\left|\frac{\partial \Psi_r(\cdot)}{\partial r}\right|\leq CI_{|x|\leq J(t)}.
$$

Finally, from (\ref{i1}) (\ref{i2}) (\ref{i3}) (\ref{i4}) (\ref{i5})  and (\ref{i1i2i3i4i5}), the claim follows. The proof is complete.

\bl(\cite[Lemma 2.1]{mps})\label{psiseco}
Suppose that $h\in C_c^2(\mR)$ is nonnegative and not identically zero. Then it holds that
\ce
\sup\limits_{x\in\mR, h(x)>0}\frac{h'(x)^2}{h(x)}\leq 2\sup\limits_{x\in\mR}|h^{\prime\prime}(x)|.
\de
\el

\section{Proof of Theorem \ref{funthe}}\label{thfunaux}

In this section, we give our verification of Theorem \ref{funthe}.

 {\bf Necessity.} On one hand, applying the It\^o formula to $V(t,\<{\bf X}_t,\varphi_t\>)$, we obtain that 
 \be
 V(t,\<{\bf X}_t,\varphi_t\>)&=&V(s,\<{\bf X}_s,\varphi_s\>)+\int_s^t\frac{\partial V(r,\<{\bf X}_r,\varphi_r\>)}{\partial r}\dif r\no\\
 &&+\int_s^t\frac{\partial V(r,\<{\bf X}_r,\varphi_r\>)}{\partial z}\left[\left<{\bf X}_r,\left(\Delta_\t\varphi_r+\frac{\partial \varphi_r}{\partial r}\right)\right>+\frac{1}{c_\t}b(r,{\bf X}_r(0))\varphi_r(0)\right]\dif r\no\\
 &&+\int_s^t\frac{\partial V(r,\<{\bf X}_r,\varphi_r\>)}{\partial z}\frac{1}{c_\t}\sigma(r,{\bf X}_r(0))\varphi_r(0)\dif B_r\no\\
 &&+\frac{1}{2}\int_s^t\frac{\partial^2 V(r,\<{\bf X}_r,\varphi_r\>)}{\partial z^2}\left(\frac{1}{c_\t}\sigma(r,{\bf X}_r(0))\varphi_r(0)\right)^2\dif r.
 \label{ito1}
 \ee
 On the other side, since $F^{\varphi}_{s,t}, V$ satisfies (\ref{defi2}), it holds that
 \be
 V(t,\<{\bf X}_t,\varphi_t\>)-V(s,\<{\bf X}_s,\varphi_s\>)=\int_s^t G_1(r, \<{\bf X}_r,\varphi_r\>)\dif r+\int_s^t G_2(r, \<{\bf X}_r,\varphi_r\>)\dif B_r.
 \label{path1}
\ee
Comparing (\ref{ito1}) and (\ref{path1}), by the uniqueness of the decomposition of semimartingales, we get 
 \ce
&&\frac{\partial V(r,\<{\bf X}_r,\varphi_r\>)}{\partial r}+\frac{1}{2}\frac{\partial^2 V(r,\<{\bf X}_r,\varphi_r\>)}{\partial z^2}\left(\frac{1}{c_\t}\sigma(r,{\bf X}_r(0))\varphi_r(0)\right)^2\\
&&+\frac{\partial V(r,\<{\bf X}_r,\varphi_r\>)}{\partial z}\left[\left<{\bf X}_r,\left(\Delta_\t\varphi_r+\frac{\partial \varphi_r}{\partial r}\right)\right>+\frac{1}{c_\t}b(r,{\bf X}_r(0))\varphi_r(0)\right]=G_1(r, \<{\bf X}_r,\varphi_r\>),\\
&&\frac{\partial V(r,\<{\bf X}_r,\varphi_r\>)}{\partial z}\frac{1}{c_\t}\sigma(r,{\bf X}_r(0))\varphi_r(0)=G_2(r, \<{\bf X}_r,\varphi_r\>).
 \de 
Furthermore, by the assumption that $\sigma(r,{\bf X}_r(0))\varphi_r(0)\not\equiv 0$ for any $r\in[0,T]$, one can get that the support of $\<{\bf X}_t,\varphi_t\>$ is $\mR$. Thus, it holds that 
\ce
&&\frac{\partial V(r,z)}{\partial r}+\frac{1}{2}\frac{\partial^2 V(r,z)}{\partial z^2}\left(\frac{1}{c_\t}\sigma(r,{\bf X}_r(0))\varphi_r(0)\right)^2\\
&&+\frac{\partial V(r,z)}{\partial z}\left[\left<{\bf X}_r,\left(\Delta_\t\varphi_r+\frac{\partial \varphi_r}{\partial r}\right)\right>+\frac{1}{c_\t}b(r,{\bf X}_r(0))\varphi_r(0)\right]=G_1(r, z),\\
&&\frac{\partial V(r,z)}{\partial z}\frac{1}{c_\t}\sigma(r,{\bf X}_r(0))\varphi_r(0)=G_2(r, z).
 \de
 Clearly, the above derived two equations are nothing but (\ref{eq2}). 

{\bf Sufficiency.} Since $V$ belongs to $C_b^{1,2}([0,T]\times\mR)$, (\ref{ito1}) is then attained. Note that $V, G_1, G_2$ satisfy the following partial differential equations
\ce\left\{\begin{array}{ll}
\partial_r V(r,z)+\frac{1}{2}\partial_z^2 V(r,z)\left(\frac{1}{c_\t}\sigma(r,{\bf X}_r(0))\varphi_r(0)\right)^2\\
+\partial_z V(r,z)\left[\<{\bf X}_r,\left(\Delta_\t\varphi_r+\frac{\partial \varphi_r}{\partial r}\right)\>+\frac{1}{c_\t}b(r,{\bf X}_r(0))\varphi_r(0)\right]=G_1(r, z),\\
\partial_z V(r,z)\frac{1}{c_\t}\sigma(r,{\bf X}_r(0))\varphi_r(0)=G_2(r, z).
\end{array}
\right.
\de
Thus, it holds that
\ce
V(t,\<{\bf X}_t,\varphi_t\>)=V(s,\<{\bf X}_s,\varphi_s\>)+\int_s^t G_1(r, \<{\bf X}_r,\varphi_r\>)\dif r+\int_s^t G_2(r, \<{\bf X}_r,\varphi_r\>)\dif B_r.
\de
That is, $F^{\varphi}_{s,t}, V$ satisfy (\ref{defi2}). The proof is complete. 

\section{Proof of Theorem \ref{pathindeth}}\label{thpathinde}

In this section, we prove Theorem \ref{pathindeth}.

Set $\psi^m(x)=p_{m^{-\frac{1}{\a}}}^\t(x)$, and $\psi^m$ weakly converges to $\d_0$ as $m\rightarrow\infty$. Thus, we replace $\varphi_t$ by $\psi^m$ in (\ref{defi2}) and obtain that
\be
\int_s^t G_1(r, \omega, \<{\bf X}_r,\psi^m\>)\dif r+\int_s^t G_2(r, \omega, \<{\bf X}_r,\psi^m\>)\dif B_r=V(t,\<{\bf X}_t,\psi^m\>)-V(s,\<{\bf X}_s,\psi^m\>).
\label{funtcond}
\ee
Besides, by Theorem \ref{funthe}, it holds that the above equality is equivalent to the following conditions
\be
&&\partial_r V(r,z)+\frac{1}{2}\partial_z^2 V(r,z)\left(\frac{1}{c_\t}\sigma(r,{\bf X}_r(0))\psi^m(0)\right)^2\no\\
&&\qquad\qquad +\partial_z V(r,z)\left[\<{\bf X}_r,\Delta_\t\psi^m\>+\frac{1}{c_\t}b(r,{\bf X}_r(0))\psi^m(0)\right]=G_1(r, z),\label{equisequ1}\\
&&\partial_z V(r,z)\frac{1}{c_\t}\sigma(r,{\bf X}_r(0))\psi^m(0)=G_2(r, z). \label{equisequ11}
\ee

Next, we observe two integrals in the left side of (\ref{funtcond}). Note that as $m\rightarrow \infty$,
$$
\<{\bf X}_r,\psi^m\>\longrightarrow {\bf X}_r(0) \quad \mbox{a.s.}.
$$
Assume that $G_1, G_2$ are uniformly bounded. Thus, by the dominated convergence theorem, it holds that
\ce
&&\int_s^t G_1(r, \omega, \<{\bf X}_r,\psi^m\>)\dif r\overset{L^1}{\longrightarrow} \int_s^t G_1(r, \omega, {\bf X}_r(0))\dif r,\\
&&\int_s^t G_2(r, \omega, \<{\bf X}_r,\psi^m\>)\dif B_r\overset{L^2}{\longrightarrow} \int_s^t G_2(r, \omega, {\bf X}_r(0))\dif B_r.
\de
Then we take the limit on two sides of (\ref{funtcond}) and attain that
\be
\int_s^t G_1(r, \omega, {\bf X}_r(0))\dif r+\int_s^t G_2(r, \omega, {\bf X}_r(0))\dif B_r=V(t,{\bf X}_t(0))-V(s,{\bf X}_s(0)).
\label{funtcond1}
\ee

Besides, as $m\rightarrow \infty$, (\ref{equisequ1}) and (\ref{equisequ11}) become 
\ce
&&\partial_r V(r,z)+\frac{1}{2}\partial_z^2 V(r,z)\left(\frac{1}{c_\t}\sigma(r,{\bf X}_r(0))\right)^2+\partial_z V(r,z)\frac{1}{c_\t}b(r,{\bf X}_r(0))=G_1(r, z),\\
&&\partial_z V(r,z)\frac{1}{c_\t}\sigma(r,{\bf X}_r(0))=G_2(r, z).
\de
Set ${\bf X}_0=x_0$, and by Subsection \ref{subspde} it holds that ${\bf X}_t(0)=X_t$ for $t\in[0,T]$ and 
\ce
&&\partial_r V(r,z)+\frac{1}{2}\partial_z^2 V(r,z)\left(\frac{1}{c_\t}\sigma(r,X_r)\right)^2+\partial_z V(r,z)\frac{1}{c_\t}b(r,X_r)=G_1(r, z),\\
&&\partial_z V(r,z)\frac{1}{c_\t}\sigma(r,X_r)=G_2(r, z).
\de

Finally, taking $g_1=G_1, g_2=G_2$ and $v=V$, we conclude that the additive functional $f_{s,t}$ is path-independent with respect to $(X_t)_{t\in[0,T]}$ if and only if both (\ref{equisequ0}) and (\ref{equisequ01}) hold. The proof is complete.

\section{An application}\label{app}

In this section, we apply Theorem \ref{pathindeth} to a class of SVEs. Since their kernels are from fractional Brownian motions, we begin with fractional Brownian motions.

First of all, let us recall some basics about fractional Brownian motions (c.f. \cite{du}). Let $B^H$ be a fractional Brownian motion with the Hurst index $H\in(0,1)$, which is a centered Gaussian process with the following covariance kernel 
$$
R_H(s,t)=\frac{V_H}{2}(s^{2H}+t^{2H}-|t-s|^{2H}),
$$
where 
$$
V_H=\frac{\Gamma(2-2H)\cos(\pi H)}{\pi H(1-2H)}.
$$
It is known that the fractional Brownian motion $B^H$ has the representation in law:
$$
B_t^H=\int_0^t K_H(t,s)\dif B_s,
$$
where $K_H(t,s)$ is the square root of the covariance operator, that is
$$
R_H(s,t)=\int_0^1 K_H(s,r)K_H(t,r)\dif r.
$$
More precisely,
$$
K_H(t,r)=\frac{(t-r)^{H-\frac{1}{2}}}{\Gamma(H+\frac{1}{2})}F(\frac{1}{2}-H, H-\frac{1}{2}, H+\frac{1}{2}, 1-\frac{t}{r})I_{[0,t)}(r),
$$
where $F$ is the Gauss hypergeometric function.

Now, we consider the following SVE:
\be
X_t=x_0+\int_0^t \tilde{K}_H(t,s)b(s,X_s)\dif s+\int_0^t \tilde{K}_H(t,s)\sigma(s,X_s)\dif B_s,
\label{fracsde}
\ee
where $\tilde{K}_H(t,s)=C(t-s)^{H-\frac{1}{2}}I_{[0,t)}(s)$, $C>0$ is a constant and $b, \sigma$ satisfy (${\bf H}^1_{b}$) (${\bf H}^1_{\sigma}$). Provided that $0<H<\frac{1}{2}$ and $\a=\frac{1}{2}-H$, by Theorem \ref{eu}, we can obtain that the equation (\ref{fracsde}) has a unique strong solution $(X_t)_{t\geq 0}$. Also assume that $\sigma\not\equiv 0$, $v$ belongs to $C_b^{1,2}([0,T]\times\mR)$ and all the derivatives of $v$ are uniformly continuous, for $\omega\in\Omega$, $g_1(\cdot,\omega,\cdot)\in C([0,T]\times\mR\mapsto\mR)$, $g_2(\cdot,\omega,\cdot)\in C([0,T]\times\mR\mapsto\mR)$ and $g_1, g_2$ are uniformly bounded. Thus, by Theorem \ref{pathindeth}, we know that the additive functional $f_{s,t}$ is path-independent with respect to $(X_t)_{t\in[0,T]}$ if and only if
 \ce
&&\partial_r v(r,z)+\frac{1}{2}\partial_z^2 v(r,z)\left(\frac{1}{c_\t}C\sigma(r,X_r)\right)^2+\partial_z v(r,z)\frac{1}{c_\t}Cb(r,X_r)=g_1(r, z),\\
&&\partial_z v(r,z)\frac{1}{c_\t}C\sigma(r,X_r)=g_2(r, z).
\de

\bigskip

\textbf{Acknowledgements:}

Two authors are very grateful to Professor Feng-Yu Wang for valuable comments. Besides, they are indebted to the referee since his suggestions and comments allowed them to improve the results and the presentation of this paper.

\end{document}